\documentclass[11pt,twoside,reqno]{amsart}
\usepackage{amsfonts}
\usepackage{}
\usepackage{latexsym,graphicx}
\usepackage[latin1]{inputenc}
\usepackage{amsmath, amsfonts, amssymb,amsthm}
\usepackage{mathrsfs}
\usepackage{color,xcolor}

\newtheorem{theorem}{Theorem}[section]
\newtheorem{definition}{Definition}[section]
\newtheorem{proposition}{Proposition}[section]
\newtheorem{lemma}{Lemma}[section]
\newtheorem{remark}{Remark}[section]

\newtheorem{corollary}[theorem]{Corollary}
\def\min{\mathop{\mathrm{min}}}

\newcommand{\N}{\mathbb{N}}

\newcommand{\R}{\mathbb{R}}

\numberwithin{equation}{section}
\allowdisplaybreaks

\begin{document}

\title[Fractional $(p,q)$-Choquard equation]{Concentrating solutions of  the fractional \\
$(p,q)$-Choquard equation with
exponential growth}

\author{Yueqiang Song$^{1},$ Xueqi Sun$^{1}$, Du\v{s}an D. Repov\v{s}$^{2,3,4}$}
\address{$^{1}$ College of Mathematics, Changchun Normal University,   Changchun, 130032,  P.R. China}
    \email{Yueqiang Song: songyq16@mails.jlu.edu.cn}
    \email{Xueqi Sun: sunxueqi1@126.com}
\address{$^{2}$ Faculty of Education, University of Ljubljana, Ljubljana, 1000, Slovenia}
\address{$^{3}$ Faculty of Mathematics and Physics, University of Ljubljana, Ljubljana, 1000, Slovenia}
\address{$^{4}$ Institute of Mathematics, Physics and Mechanics, Ljubljana, 1000, Slovenia}
\email{Du\v{s}an D. Repov\v{s}: dusan.repovs@guest.arnes.si}

\thanks{$^*$ Corresponding author: Du\v{s}an D. Repov\v{s}}
\thanks{{\it Mathematics Subject Classification (2020):} 35A15, 35A23, 35J35, 35J60, 35R11.}
\thanks{{\it Keywords}: Fractional double phase operator, Critical exponential growth,  Mountain Pass Theorem,
    Trudinger-Moser inequality, Variational method.}

\begin{abstract}
This article deals with the following fractional $(p,q)$-Choquard
equation with exponential growth of the form:
    $$ \varepsilon^{ps}(-\Delta)_{p}^{s}u+\varepsilon^{qs}(-\Delta)_q^su+ Z(x)(|u|^{p-2}u+|u|^{q-2}u)=\varepsilon^{\mu-N}[|x|^{-\mu}*F(u)]f(u) \ \ \mbox{in} \ \ \mathbb{R}^N,$$
    where $s\in (0,1),$ $\varepsilon>0$ is a parameter,  $2\leq p=\frac{N}{s}<q,$ and $0<\mu<N.$
    The nonlinear function $f$ has an exponential growth at infinity and the continuous potential
    function $Z$ satisfies suitable
    natural conditions. With the help of the Ljusternik-Schnirelmann category theory and variational
    methods, the multiplicity and concentration of positive solutions  are obtained for  $\varepsilon>0$ small
enough. In a certain sense, we generalize some  previously known results.
\end{abstract}

\baselineskip=13truept
\maketitle
\pagestyle{myheadings}
\markboth{}{}

\section{Introduction}\label{s1}
\def\theequation{1.\arabic{equation}}
\setcounter{equation}{0}

In this paper, we consider the multiplicity and concentration  of
solutions for the following fractional $(p,q)$-Choquard problem in
$\mathbb{R}^N$:
\begin{align}\label{pr2}
\varepsilon^{ps}(-\Delta)_{p}^{s}u+\varepsilon^{qs}(-\Delta)_q^su+Z(x)(|u|^{p-2}u+|u|^{q-2}u)
=\varepsilon^{\mu-N}[|x|^{-\mu}*F(u)]f(u)\tag{$\mathcal Q$},
\end{align}
where $\varepsilon$ is small positive parameter, $0<\mu<N$, $0<s<1$,
$2\le p=\frac{N}{s}<q,$ the continuous potential $Z$ is bounded from
below by $Z_0>0$, the nonlinearity  $f$ has an exponential critical
growth at infinity, and $(-\Delta)_{\wp}^{s}$ $(\wp\in \{p,q\})$ is the
fractional $\wp$-Laplace operator defined by
$$ (-\Delta)_{\wp}^{s}u(x)=2\lim_{r\to 0^{+}}\int\limits_{\mathbb R^N\setminus B_{r}(x)}
\frac{|u(x)-u(y)|^{\wp-2}(u(x)-u(y))}{|x-y|^{N+\wp s}}dy \ 
\hbox{ for every } x\in \mathbb R^{N},
$$
 up to a normalization
constant in the integral, where $u \in C_{0}^{\infty}(\mathbb R^{N})$
 and $B_{r}(x)$ denotes the ball with center $x$ of
  radius $r>0.$

  Many scholars have studied
fractional and nonlocal operators because of their applications in various contexts,
for example, in optimization, finance, crystal dislocations, phase transitions, etc. For more on
these topics, we refer to
Ambrosio \cite{Amb1}
and
di Nezza et al. \cite{Ne1}.

We shall assume that the potential function $Z$ and the nonlinearity $f$ satisfy the following conditions:
\begin{itemize}
\item[$(\mathcal {Z}_{1})$]  There exists $Z_{0}>0$ such that $Z(x) \geq {Z_0},$ for every ${x \in {\mathbb{R}^N}};$

\item[$(\mathcal {Z}_{2})$]  There exists an open bounded set $\Omega\subset\mathbb{R}^{N}$ such that
\[{Z_0} = \mathop {\inf }\limits_{x \in \Omega } Z\left( x \right)<\min\limits_{x \in \partial \Omega } Z\left( x \right).\]

\item [$(f_1)$]
  $f$ is a continuous function
such that $f(t)=0$
and
 for every $t\le 0$
and
 every $q_1$, $q_2$, such that
$$q_1\ge q,\qquad q_2\ge \frac{N}{s},$$
there exist  real  numbers $a_1>0$, $a_2>0,$
and $\beta_0$, with $0<\beta_0<\alpha_{*}(s,N)$, such that
$$ f'(t)\le a_1|t|^{q_1-2}+a_2\mathcal {H}_{N,s}(\beta_0|t|^{N/(N-s)})|t|^{q_2-2} \
{\mbox{ for every }t\ge0,}$$
where (see
Parini and Ruf \cite{PR18}
and Zhang \cite{CZ})
$$\mathcal {H}_{N,s}(t)=e^{t}-\sum_{j=0}^{j_p-2}\frac{t^{j}}{j!},\quad
 j_p=\min\{j\in \mathbb N: j\ge p\},\mbox{ and }$$
 \begin{equation}\label{alpha*}
 \begin{gathered}\alpha_{*}=
 \alpha_{*}(s,N)=N\left[\frac{2(N\omega_N)^2\Gamma(p+1)}{N!}
 \sum_{k=0}^{\infty}\frac{(N+k-1)!}{k!}\frac1{(N+2k)^p}\right]^{s/(N-s)},\\
 \quad \omega_N=\frac{\pi^{N/2}}{\Gamma(1+N/2)}.\end{gathered}
 \end{equation}

 \item [$(f_2)$]
 $\lim\limits_{t\rightarrow 0^+}\frac{f'(t)}{t^{q-2}}=0.$

  \item [$(f_3)$]
  There exists $\theta>q$ such that $f(t)t\ge \theta F(t)>0,$ for every $t>0,$ where $F(t)=\int_{0}^{t}f(\tau)d\tau.$

 \item [$(f_4)$]
 There exists $\gamma_1>0$ large enough such that $F(t)\ge \gamma_1 |t|^{\theta},$ for every $t\ge 0,$ where $\theta$ is as  given in $(f_3)$.

 \item [$(f_5)$]
  The function $t\mapsto
f(t)t^{1-q}$ is strictly increasing on
$\mathbb R^+=(0, +\infty).$
\end{itemize}

Once $s=1$, problem \eqref{pr2} reduces to a typical $(p,q)$-elliptic equation:
\begin{align}\label{m2c}
    -\Delta_pu- \Delta_qu+Z(\varepsilon x)(|u|^{p-2}u+|u|^{q-2}u)=H(u) \text{ in }\; \mathbb R^N,
\end{align}
where $H$ is nonlinear reaction, $\Delta_{\wp}u=\text{div}(|\nabla
u|^{\wp-2}\nabla u),$ and $\wp\in \{p,q\}.$
It has been widely studied in physics, biophysics, plasma physics, chemical reaction design and elsewhere.
For more physical examples, we refer
to
Antontsev and Shmarev \cite{AS},
Benci et al. \cite{BDaFP},
and
Cherfil and Il'yasov \cite{CV},
and the references therein.
The multiple phases equation was proposed in the study of the Born-Infeld equation
(see Bonheure et al. \cite{BAP},
Born and Infeld \cite{BL},
and
Br\'ezis and Lieb \cite{BL1}),
which
 models  electromagnetic fields, electrostatics and electrodynamics, and was a model based on the Maxwell-Lagrangian density
$$ -\text{div}\left(\dfrac{\nabla u}{(1-2|\nabla u|^2)^{1/2}}\right)=h(u)\;\text{in}\;\mathbb R^N.$$

When $p=q$, problem \eqref{pr2} becomes the fractional $p$-Laplace Choquard equation of the form:
\begin{align}\label{eqm2}
    \varepsilon^{ps}(-\Delta)_p^su+Z(x)|u|^{p-2}u=\varepsilon^{\mu-N}
    \Big[|x|^{-\mu}*F(u)\Big]f(u(x))\;\;\text{in }\;\mathbb R^N,
\end{align}
where $\varepsilon>0$ is a sufficiently small parameter, typically the Planck
constant and $F(t)=\int_0^t f(\tau)d \tau$. We say that a solution of problem \eqref{eqm2}
is {\em semi-classical} if
 $\varepsilon\to 0^{+}$. From the physics point of view, the semi-classical solution is also a solutions of problem \eqref{eqm2}, when  $\varepsilon\to
{0^+}$.
Floer and Weinstein~\cite{FW86}
established the existence of semi-classical solutions of problem \eqref{eqm2}.
 A special  form of problem \eqref{eqm2} is
\begin{align}\label{pr1}
    -\varepsilon^2\Delta
    u+Z(x)u=\varepsilon^{\mu-N}[|x|^{-\mu}*F(u)]f(u)\;\
    \text{in}\;\mathbb R^N \
    \hbox{ for every }
     0<\mu<N,
\end{align}
where $f$ is a nonlinear reaction. It is worth noting that problem \eqref{pr1} was introduced
in the theory of the Bose-Einstein condensation and used to describe the
finite-range many body interactions between particles. There are already many works on this topic.
 By  variational methods, Alves et al.~\cite{ACTY} considered the concentration
solutions of problem \eqref{pr1} in $\mathbb{R}^2,$ where $f$ has exponential critical
growth and $Z$ satisfies some appropriate conditions.  Once $F(u)=|u|^{p}$ in problem \eqref{pr1},
one obtains the following Choquard equation
\begin{align}\label{ctma0}
    -\Delta u+Z(x)u=(I_{\alpha}*|u|^{p})|u|^{p-2}u\;\;\text{in}\;\mathbb R^N,
\end{align}
where $I_{\alpha}$ is the Riesz potential,
$\Gamma$ is the Gamma function, and $Z$ is a potential function.

 When
$p=\alpha=2$, $N=3,$ and $Z(x)=\nu,$ problem \eqref{ctma0}
reduces to the Choquard-Pekar type equation
\begin{align}\label{ctma1}
    -\Delta u+\nu u=(I_2*u^2)u, \hbox{ for every } x\in\mathbb R^3,
\end{align}
which was
proposed in 1976 by  Lieb \cite{L76},  in order to describe
an electron trapped in its own hole. Problem \eqref{ctma1} is called the Schr\"{o}dinger-Newton equation. Inspired by the work  of Lieb~\cite{L76} and Lions
\cite{PL80}, many researchers have studied the Choquard equation by variational
methods.

Recently, these methods have become more useful for establishing
the existence of weak solutions of the Choquard equations. 
For example, Chen and Yang \cite{chen1} studied the following Choquard equation with upper critical exponent on a bounded domain
$$ -\Delta u = \mu f(x)|u|^{p-2}u + g(x)(I_\alpha^\ast(g|u|^{2_\alpha^\ast}))|u|^{2_\alpha^\ast-2}u\ \hbox{ for every }  x \in \Omega,$$
where $\mu > 0$ is a parameter, $N > 4$, $0 < \alpha < N$, $I_\alpha$ is the Riesz potential, $\frac{N}{N-2}<p<2$, $\Omega$ is a bounded
domain with smooth boundary, and $f$ and $g$ are continuous functions. For $\mu$ small enough, with the help of variational
methods, they established the relationship between the number of solutions and the profile of potential $g$.

Yang and Zhao~\cite{Yang} studied the singularly perturbed fractional Choquard equation
\begin{equation}\label{3a}
\varepsilon^{2s}(-\Delta)^su+V(x)u=\varepsilon^{\mu-3}\left(\int_{\mathbb{R}^3}
\frac{|u(y)|^{2_{\mu,s}^\ast}+F(u(y))}{|x-y|^{\mu}}dy\right)(|u|^{2_{\mu,s}^\ast-2}u+\frac{1}{2_{\mu,s}^\ast}f(u))
\hbox{ in } \mathbb{R}^N,
\end{equation}
 where $2_{\mu,s}^\ast=\frac{6-\mu}{3-2s}$ is the critical exponent in the sense of the
Hardy-Littlewood-Sobolev inequality and the continuous function $f$ satisfies subcritical growth conditions. By variational
methods, penalization techniques and the Lyusternik-Schnirelmann theory, the authors established the
multiplicity and concentration behaviour of solutions  for problem \eqref{3a}.

We need to point out some recent results: Zuo et al. \cite{zuo} developed a variational approach, based on the scaling function method to solve optimization problems. Here, the authors dealt with the mass subcritical case, and referred to the fractional framework setting. Zhang et al. \cite{zl} considered a class of fractional parabolic equation with general nonlinearities. The authors established monotone increasing property of the positive solutions in one direction. Based on this, nonexistence of the solutions was demonstrated, via a contradiction argument.
For more information, we refer to
Chen et al. \cite{CLY,CTS22},
Cingolani and Tanaka \cite{Cin},
Clemente et al. \cite{CAB},
and
B\"{o}er and Miyagaki \cite{EM21},
 and  the references therein.

For  fractional $(p,q)$-Laplace problems, some interesting existence and multiplicity results have emerged in recent years.
Zhang et al.~\cite{ZZR23} studied multiplicity and
concentration solution for the double phase equation  in $\mathbb
R^N,$ especially, they assumed the nonlinearity of $f\in C^1(\mathbb{R}^N)$ and that the continuous potential
function  satisfies the global condition.  Later, using
the penalization method, the Ljusternik-Schnirelmann theory,
and variational  methods,
Ambrosio~\cite{A0} first studied existence of
multiple solutions and concentration of the $(p,q)$-fractional
Choquard equation
\begin{align}\label{m1a}
    (-\Delta)_p^su+(-\Delta)_q^{s}u+ Z(\varepsilon x)(|u|^{p-2}u+|u|^{q-2}u)=(|x|^{-\mu}*F(u))f(u)
    \ \ \mbox{in} \ \ \mathbb{R}^N,
\end{align}
where $\mu\in [0,ps)$, $0<s<1,$ $1<p<q<\frac{N}{s},$  $f$ has subcritical growth and the potential function $Z$ satisfies the local conditions.
Molica
Bisci et al. ~\cite{BTR23} extended the results of Zhang
et al.~\cite{ZZR23} to the fractional Choquard  problem
\eqref{m1a}, where the potential function $Z(x)$ satisfies the global condition.
With the help of  the Ljusternik-Schnirelmann category theory and variational methods,
Liang et al.~\cite{Liang1} explored the multiplicity and concentration behaviors of solutions for  the $(p,q)$ fractional Choquard equation with exponential growth.

To the best of our knowledge, there are  no known results concerning problem
\eqref{pr2}, when the continuous function $f$ has the  exponential growth behavior at infinity in
the sense of Trudinger-Moser. Inspired by the results of Liang et al.~\cite{Liang1}, 
 we show in this paper the existence and concentration behavior of
solutions of problem \eqref{pr2} involving exponential growth.
In comparison to Liang et al.~\cite{Liang1},
 we assume the nonlinearity and function  $f$ is supposed to be only continuous, which makes the
corresponding Nehari manifold possibly nondifferentiable. Therefore, we cannot directly  use the differentiability of the Nehari manifold.
Furthermore, it is not possible to apply  the
Ljusternik-Schnirelmann category theory on the Nehari manifold in
order to obtain the multiplicity of solutions for problem
\eqref{pr2}. Whereas  Liang et al.~\cite{Liang1} have $f\in C^1(\mathbb{R}),$
 we need to apply some other techniques to overcome this
difficulty.

In addition, this is the first time that  problem
\eqref{pr2} with the Trudinger-Moser nonlinearities has been studied in both cases: $s\in(0,1)$
 and $s\rightarrow1^-$. Comparing with Ambrosio~\cite{A0},
he studied the subcritical growth $p<\frac{N}{s}$
and  the local case. We obtain the Sobolev embedding from
$W^{s,p}(\mathbb R^N)$ into $L^t(\mathbb R^N),$ for every $t\in
[p,p_s^{*}]$. However,  in this paper we consider the case $N=ps$, so
the embedding $W^{s,p}(\mathbb R^N)\hookrightarrow
L^{\infty}(\mathbb R^N)$ may not exist.
To overcome this obstacle, it is essential to apply the fractional Trudinger-Moser inequality,
which
 is the main difference with
Ambrosio  \cite{A0},
Molica Bisci et al. \cite{BTR23},
and
Zhang et al. \cite{ZZR23}.
The other major challege which we encountered,  is the loss of compactness of the
Palais-Smale sequences associated with the underlying functionals,
corresponding  to problems \eqref{pr2} and \eqref{pr3}, so we use
certain analytical techniques to overcome this obstacle.

\begin{definition}\label{de1}
We denote the category of a set $A$ with
respect to a set $B$  by $\text{cat}_B(A)$ as the least integer $k$ such that
$A\subset A_1\cup \cdots\cup A_k,$ where each $A_i$, $i=1,\cdots,k$, is a
closed and contractible subset of $B.$ We set
$\text{cat}_{B}(\emptyset)=0$ and $\text{cat}_B(A)=\infty$ if there
is no integer with the above property.
\end{definition}
  Let
\begin{align}\label{M}\mathscr{M}=\{x\in\mathbb R^N\,:\, Z(x)=Z_0\}\end{align}
and for every $\delta>0$ define
$$ \mathscr{M}_{\delta}=\{x\in\mathbb R^N\,:\, \mbox{dist}(x,\mathscr{M})\le \delta\}.$$
We are now ready to state the main results of this paper.

\begin{theorem}\label{th2}
Suppose  that  conditions $(\mathcal {Z}_1),$ $(\mathcal {Z}_2),$ and $(f_1)-(f_5)$ are satisfied. Then for every $\delta>0,$ there
exists $\varepsilon_{\delta}>0$ such that problem \eqref{pr2} has
at least $\textrm{cat}_{\mathscr{M}_{\delta}}(\mathscr{M})$ positive (weak) solutions
for every  $\varepsilon>0$ satisfying $\varepsilon<\varepsilon_{\delta}.$
Furthermore, let $w_{\varepsilon}$ be a solution of problem \eqref{pr2} and
$\zeta_{\varepsilon}$  its global maximum.  Then, up to a
subsequence,
    $\zeta_{\varepsilon}\to y\in \mathscr{M}$ and
    $ \lim\limits_{\varepsilon\to 0^{+}} Z(\zeta_{\varepsilon})=Z_{0}.$
\end{theorem}

\begin{theorem}\label{th2a}
Suppose  that  conditions $(\mathcal {Z}_1),$ $(\mathcal {Z}_2),$ and $(f_1)-(f_5)$ are satisfied.
Let
$w_{\varepsilon}$ be a solution of problem \eqref{pr2}, which
exists by Theorem~$\ref{th2}$, and let $\zeta_{\varepsilon}$ be its
global maximum. Then $u_{\varepsilon}(x)=w_{\varepsilon}(\varepsilon
x+\zeta_{\varepsilon})$ converges strongly in $W^{s,p}(\mathbb
R^N)\cap W^{s,q}(\mathbb R^N)$ to a ground state solution $u$ of
 the following problem
$$(-\Delta)_p^su+(-\Delta)_q^su+Z_{0}(|u|^{p-2}u+|u|^{q-2}u)=[|x|^{-\mu}*F(u)]
f(u)\mbox{ in }\mathbb R^N.$$
\end{theorem}

\begin{remark}\label{rem2}
{\em As $s\to 1^{-1}$, condition $(f_1)$ reduces to  the following condition}  \\

\noindent $(f_1)'$ {\em  Continuous function $f(t)$ vanishes for every $t\in(-\infty,0)$}
{\em and for every $q_1$, $q_2$, with
$$q_1\ge q,\qquad q_2>N,$$
 there exist   constants $a_1>0$, $a_2>0$
and $\alpha_0$, with $0<\beta_0<\alpha_{*}$, such that
$$ f(t)\le a_1|t|^{q_1-1}+a_2\mathcal {H}_{N}(\beta_0|t|^{N/(N-1)})|t|^{q_2-1}\
{\hbox{ for every } t\in\mathbb R^+_0},$$
where
$$\mathcal {H}_{N}(t)=e^{t}-\sum_{j=0}^{N-2}\dfrac{t^{j}}{j!},\quad
\quad 0< \alpha_{*}\le \alpha_{*}(1,N)=\lim_{s\to 1^{-}}
\alpha_{*}(s,N) $$ and $\alpha_{*}(s,N)$ is given
in~\eqref{alpha*}.}
\end{remark}

Invoking Theorem \ref{th2} and Theorem \ref{th2a} for the case when  $s\to 1^{-1},$
we get the following
results for problem \eqref{pr2}, respectively.

\begin{corollary}\label{cornew10}
Suppose that conditions $(\mathcal {Z}_1), (\mathcal {Z}_2),$
 $(f_1)',$ and  $(f_2)-(f_5)$ hold. Then for every $\delta>0,$
there exists $\varepsilon_{\delta}>0$ such that problem
\eqref{pr2} has at least $\textrm{cat}_{\mathscr{M}_{\delta}}(\mathscr{M})$
positive (weak) solutions for every $\varepsilon \in (0,\varepsilon_{\delta}).$
Furthermore, let $w_{\varepsilon}$ be a solution
of problem \eqref{pr2}
and $\zeta_{\varepsilon}$  its
global maximum.
Then, up to a subsequence,
$\zeta_{\varepsilon}\to y\in \mathscr{M}$ and $ \lim_{\varepsilon\to 0^{+}} Z(\zeta_{\varepsilon})=Z_{0}.$
\end{corollary}

\begin{corollary}\label{cornew1}
 Suppose that  conditions $(\mathcal {Z}_1), (\mathcal {Z}_2),$ $(f_1)'$, and $(f_2)-(f_5)$ hold.
 If $w_{\varepsilon}$ is a solution of problem \eqref{pr2}, which
exists by   Corollary~$\ref{cornew10}$, and $\zeta_{\varepsilon}$ is
its global maximum, then
$u_{\varepsilon}(x)=w_{\varepsilon}(\varepsilon
x+\zeta_{\varepsilon})$ converges strongly in $W^{1,N}(\mathbb
R^N)\cap W^{1,q}(\mathbb R^N)$ to a ground state solution $u$ of
$$-\Delta_Nu-\Delta_qu+Z_{0}(|u|^{N-2}u+|u|^{q-2}u)=[|x|^{-\mu}*F(u)]f(u)
\text{ in }\;\mathbb R^N$$
 and there exist $c>0, C>0$ such that $|w_{\varepsilon}(x)|\le Ce^{-c|x-\zeta_{\varepsilon}|/\varepsilon},$ for every $x\in\mathbb R^N.$
\end{corollary}

The organization of this paper is as follows.
 In Section~\ref{s2}, we introduce some  notations and recall certain technical
results which will be needed in the paper.
 In Section~\ref{s3}, we study the autonomous problem
\eqref{ct1a} associated to problem \eqref{pr2}.
In Section~\ref{s4}, we deal with the auxiliary problem \eqref{m4}. In addition, we also verify that the Palais-Smale
condition holds for its  energy functional and apply some new tools  to obtain a multiplicity result.
 In Section~\ref{s5}, we establish the multiplicity of solutions for the modified problem and complete the proof of the main results.

\section{Preliminaries}\label{s2}
\def\theequation{2.\arabic{equation}}
\setcounter{equation}{0}

In this section, we state some results and notions which will be used later.
For all other background material we refer to the comprehensive monograph by Papageorgiou et al.  \cite{PRR}.

Let us recall that
$p=\frac{N}{s}$ in problem \eqref{pr2}. For $1<p<\infty,$ we define the fractional
Sobolev space $W^{s, p}(\mathbb R^N)$ as
$$ W^{s, p}(\mathbb R^N)=\Big\{u\in L^{p}(\mathbb R^N)\,:\, [u]_{s, p}=\left(\iint_{\mathbb R^{2N}}\dfrac{|u(x)-u(y)|^{p}}{|x-y|^{2N}}dxdy\right)^{1/p}<\infty\Big\} $$
and we endow it  with the norm
$$ \|u\|_{W^{s,p}(\mathbb{R}^N)}=\Big(\|u\|_{L^p(\mathbb{R}^N)}^{p}+[u]_{s, p}^{p}\Big)^{1/p}.$$
By Pucci et al.~\cite[Lemma 10]{PXZ}, the fractional
Sobolev space $W^{s, p}(\mathbb R^N)$ is a uniformly
convex Banach space.

Now, fix $\vartheta>0$ and endow $W^{s, p}(\mathbb R^N)$ with the
norm
$$ \|u\|_{\vartheta,W^{s,p}(\mathbb R^N)}=\Big(\vartheta\|u\|_{L^p(\mathbb{R}^N)}^{p}+[u]_{s, p}^{p}\Big)^{1/p}.$$
Obviously, the norm $\|\cdot\|_{W^{s,p}(\mathbb{R}^N)}$ is equivalent to
$\|\cdot\|_{\vartheta,W^{s,p}(\mathbb R^N)}$ on
$W^{s, p}(\mathbb R^N).$

Let conditions $(\mathcal {Z}_1)$ and $(\mathcal {Z}_2)$ be satisfied.  We denote by
$W_{Z,\varepsilon}^{s,p}(\mathbb R^N)$ the completion of
$C_0^{\infty}(\mathbb R^N),$ with norm
\begin{align*}
    \|u\|_{W_{Z,\varepsilon}^{s,p}(\mathbb R^N)}=\Big([u]_{s,p}^{p}+ \|u\|_{p,Z,\varepsilon}^{p}\Big)^{1/p}, \quad \|u\|_{p,Z,\varepsilon}^{p}=\int\limits_{\mathbb R^N}
    Z(\varepsilon x)|u(x)|^{p}dx\
    \hbox{ for every }
    \varepsilon>0.
\end{align*}

Pucci et al. ~\cite[Lemma 10]{PXZ}
showed
 that $W_{Z,\varepsilon}^{s,p}(\mathbb R^N)$
 is also a uniformly convex Banach space for $1<p<\infty$. Moreover,
$W_{Z,\varepsilon}^{s,p}(\mathbb R^N)$ is a reflexive Banach space. Invoking
conditions $(\mathcal {Z}_1)-(\mathcal {Z}_2),$ and Di Nezza et al.~\cite[Theorem~6.9]{Ne1}, we
obtain the continuous embedding $W_{Z,\varepsilon}^{s,p}(\mathbb
R^N)\hookrightarrow L^{\nu}(\mathbb R^N)$  for arbitrary
$\nu\in[N/s, \infty).$

Let $s\in(0,1)$ and $p,q\in(1,\infty).$
The natural solution space of problem  \eqref{pr2} is defined as
$$\mathcal {W}_{\varepsilon}=W_{Z,\varepsilon}^{s,p}(\mathbb R^N)\cap
W_{Z,\varepsilon}^{s,q}(\mathbb R^N) $$ and is equipped with the norm
$$ \|u\|_{\mathcal {W}_{\varepsilon}}=\|u\|_{W_{Z,\varepsilon}^{s,p}(\mathbb R^N)}+ \|u\|_{W_{Z,\varepsilon}^{s,q}(\mathbb R^N)}.$$
With  the aid of above definitions, assumptions $(\mathcal {Z}_1)-\mathcal {Z}_2),$ and the fact that $p=\frac{N}{s}$, it is easy to get the continuous embeddings
$$\mathcal {W}_{\varepsilon}\hookrightarrow W_{Z,\varepsilon}^{s,p}(\mathbb R^N)\hookrightarrow W^{s,p}(\mathbb R^N)\hookrightarrow L^{\nu}(\mathbb R^N)\
 \text{ for every }  \nu\in[{N}/{s}, \infty).$$ Therefore, there exist the best constants
$$S_{\nu,\varepsilon}
=\mathop{\inf}\limits_{\substack{u\in \mathcal {W}_{\varepsilon}\\u\neq 0}}
\frac{\|u\|_{\mathcal {W}_{\varepsilon}}}{\|u\|_{L^{\nu}(\mathbb R^N)}}\
\hbox{ for every }
\nu\in [{N}/{s}, \infty).
$$

By a change of variable $x\mapsto \varepsilon x,$  problem
\eqref{pr2} becomes equivalent to the following equation
\begin{equation}\label{pr3}
    (-\Delta)_{N/s}^{s}u+(-\Delta)_q^su+Z(\varepsilon
    x)(|u|^{\frac{N}{s}-2}u+|u|^{q-2}u)=[|x|^{-\mu}*F(u)]f(u)
   \tag{$\mathcal {Q}_{\varepsilon}$}\
   \hbox{ in }
   \mathbb R^N
\end{equation}
 which is variational and the (weak) solutions of problem \eqref{pr3} satisfy
the  following definition.
\begin{definition}\label{dn1}
{\rm  Let $u\in \mathcal {W}_{\varepsilon}.$ If for every $\varphi\in \mathcal {W}_{\varepsilon},$ we have
    \begin{align*}
        \sum_{\wp\in \{p,q\}}&\iint\limits_{\mathbb R^{2N}}\frac{|u(x)-u(y)|^{\wp-2}(u(x)-u(y))(\varphi(x)-\varphi(y))}{|x-y|^{N+\wp s}}dxdy\\
        &+\int\limits_{\mathbb R^N}Z(\varepsilon x)(|u|^{p-2}u+|u|^{q-2}u)\varphi dx=\int\limits_{\mathbb R^N}\int\limits_{\mathbb R^N}\dfrac{F(u(y))}{|x-y|^{\mu}}f(u(x))\varphi(x)dxdy,
    \end{align*}
\text{ then $u$ is  called a weak {\em solution} of~problem \eqref{pr3}}.
    }
\end{definition}

\section{The autonomous problem $(\mathcal {Q}_{\vartheta})$}\label{s3}
\def\theequation{3.\arabic{equation}}
\setcounter{equation}{0}

Fix $\vartheta>0$. In this section, we shall consider the
autonomous problem $(\mathcal {Q}_{\vartheta})$, associated  with problem
\eqref{pr2}, that is
\begin{align*}\label{ct1a}
    (-\Delta)_{N/s}^{s}u+(-\Delta)_q^su+\vartheta
    (|u|^{\frac{N}{s}-2}u+|u|^{q-2}u)=[|x|^{-\mu}*F(u)]f(u)\;\ \text{in}\;
    \mathbb R^N.\tag{$\mathcal {Q}_{\vartheta}$}
\end{align*}
 We consider the Euler-Lagrange functional  $\mathcal {E}_{\vartheta}: W^{s,N/s}(\mathbb
R^N)\cap W^{s,q}( \mathbb R^N)\to \mathbb R$ corresponding to problem
$(\mathcal {Q}_{\vartheta})$ as follows
\begin{align}\label{Jeta}
 \mathcal E_{\vartheta}(u)&=\frac{1}{p}\|u\|_{\vartheta,W^{s,p}(\mathbb R^N)}^{p}+\dfrac{1}{q}\|u\|_{\vartheta,W^{s,q}(\mathbb R^N)}^{q}-\frac{1}{2}\int\limits_{\mathbb R^N}\mathcal {L}_\mu(x)F(u(x))dx,
\end{align}
where
$$\mathcal {L}_\mu(u)(x)=\int\limits_{\mathbb
R^N}\frac{F(u(y))}{|x-y|^{\mu}}dy.$$ Here, $\mathcal {W}=W^{s,N/s}(\mathbb
R^N)\cap W^{s,q}( \mathbb R^N)$ is the Banach space, with the norm
$$\|u\|=\|u\|_{W^{s,p}(\mathbb R^N)}+\|u\|_{W^{s,q}(\mathbb R^N)}.$$
We also endow $\mathcal {W}$ with  the equivalent norm
$$\|u\|_{\vartheta,\mathcal {W}}=\|u\|_{\vartheta, W^{s,p}(\mathbb R^N)}+\|u\|_{\vartheta,W^{s,q}(\mathbb R^N)}.$$
Thus, $\mathcal {W}$ is a uniformly convex  Banach space
and so $\mathcal {W}$ is also a
reflexive Banach space. By  Di Nezza et al.~\cite[Theorem~6.9 ]{Ne1}, we obtain
the continuous embeddings
$$\mathcal {W}\hookrightarrow W^{s,N/s}(\mathbb R^N)\hookrightarrow
L^{\nu}(\mathbb R^N)\ \hbox{ for every }  \nu\in[{N}/{s},\infty).$$
 Hence,
there exists the best constant $A_{\nu,\eta}>0$  given by
$$A_{\nu,\vartheta}=\inf\limits_{\substack{u\in \mathcal {W}\\u\neq 0}}
\frac{\|u\|_{\vartheta,\mathcal {W}}}{\|u\|_{L^{\nu}(\mathbb R^N)}}\
\hbox{ for every }
\nu\in [{N}/{s}, \infty).$$

\begin{lemma}(see Zhang \cite[Theorem 1.1]{CZ})\label{lm1}
Let $s\in (0,1)$ and $N=sp.$ Then for every $\alpha$, with
$0<\alpha<\alpha_{*}\le \alpha_{*}(s,N),$
\begin{gather*}\sup\limits_{\substack{v\in W^{s,p}(\mathbb R^N)\\||v||_{W^{s,p}(\mathbb R^N)}\le 1}}\,\int\limits_{\mathbb R^N}\mathcal {H}_{N,s}(\alpha |v|^{N/(N-s)})dx<\infty,\\
\mbox{where\quad}\mathcal {H}_{N,s}(t)=e^{t}-\sum_{j=0}^{j_p-2}\dfrac{t^j}{j!},\qquad
j_p=\min \{j\in \mathbb N\,:\, j\ge p\}.\end{gather*}
Moreover, for $\alpha>\alpha_{*}(s,N),$
    $$ \sup\limits_{\substack{v\in W^{s,p}(\mathbb R^N)\\||v||_{W^{s,p}(\mathbb R^N)}\le 1}}\,
    \int\limits_{\mathbb R^N}\mathcal {H}_{N,s}(\alpha |v|^{N/(N-s)})dx=\infty,$$
    where $$ \alpha_{*}(s,N)=N\left(\dfrac{2(N\omega_N)^2\Gamma(p+1)}{N!}\sum_{k=0}^{+\infty}
    \dfrac{(N+k-1)!}{k!}\dfrac{1}{(N+2k)^{p}}\right)^{s/(N-s)}=N(\gamma_{s,N})^{s/(N-s)}.$$
\end{lemma}

\begin{remark}\label{rm1a}
 {\rm   In Lemma \ref{lm1}, if we take the norm $\|.\|_{\vartheta,W^{s,p}(\mathbb R^N)}$ in  $W^{s,N/s}(\mathbb R^N),$ then
    $$ (\max\{1,\vartheta\})^{-1/p}\|v\|_{\vartheta,W^{s,p}(\mathbb R^N)}\le \|v\|_{W^{s,p}(\mathbb R^N)}\le (\min\{1,\vartheta\})^{-1/p}||v||_{\vartheta,W^{s,p}(\mathbb R^N)},    
    v\in W^{s,N/s}(\mathbb R^N).
    $$
Moreover,
$$ \sup\limits_{\substack{v\in W^{s,p}(\mathbb R^N)\\\|v\|_{\vartheta,W^{s,p}
(\mathbb R^N)}\le (\min\{1,\vartheta\})^{s/N}}}\,
\int\limits_{\mathbb R^N}\mathcal {H}_{N,s}(\alpha |v|^{N/(N-s)})dx<\infty\
\hbox{ for every }
\alpha,
0< \alpha<\alpha_{*}\le {\alpha_{*}(s,N)}.
$$}
\end{remark}

\begin{lemma}[The Hardy-Littlewood-Sobolev inequality, see Lieb \cite{L76, Lie}]\label{lm2n}
Let $r$, $t>1,$ and $0<\mu<N$ such that
$$\dfrac{1}{r}+\dfrac{\mu}{N}+\dfrac{1}{t}=2.$$
Then there exists a sharp constant $C(r,N,\mu,t)>0$ such that
$$\int\limits_{\mathbb R^N}\int\limits_{\mathbb R^N}\dfrac{{|u(x)v(y)|}}{|x-y|^{\mu}}dxdy\le C(r,N,\mu,t)\|u\|_{L^r(\mathbb R^N)}\|v\|_{L^{t}(\mathbb R^N)}\
\hbox{ for every }
u\in L^{r}(\mathbb R^N),
v\in L^{t}(\mathbb R^N).$$
\end{lemma}

Since $r=t$, we can use Lemma~\ref{lm2n} and the following equality
$$\frac{2}{t}+\frac{\mu}{N}=2,$$
that is, when
$ t=\frac{2N}{2N-\mu},$
then the integral
$$\int\limits_{\mathbb R^N}[|x|^{-\mu}*F(u)]F(u)dx \ \  \mbox{for every}   \ F(u)=|u|^{q}$$
is well-defined on $L^{t}(\mathbb R^N)$, with $t=2N/(2N-\mu),$ along
every $u\in W^{s,N/s}(\mathbb R^N)$, provided that $qt\ge
\frac{N}{s},$ due to the  continuous embedding $W^{s,N/s}(\mathbb
R^N)\hookrightarrow L^{\nu}(\mathbb R^N),$ for every $\nu\in
[N/s,\infty)$. Hence,
$$q\ge \frac{N}{st}=\frac{2N-\mu}{2s}.$$
Consequently, Lemma~\ref{lm1}, assumption $(f_1)$
and the fact that $C_0^{\infty}(\mathbb R^N)$ is dense
in $W^{s,p}(\mathbb R^N),$ imply that~$\mathcal {E}_{\vartheta}$ is well-defined on $\mathcal {W}$
and of class $C^1(\mathcal {W})$.
Furthermore, for every $u\in \mathcal {W},$
\begin{align*}
    \langle \mathcal {E}'_{\vartheta}(u),\varphi\rangle&=\sum_{\wp\in \{p,q\}}\int\limits_{\mathbb R^{2N}}\frac{|u(x)-u(y)|^{\wp-2}(u(x)-u(y))(\varphi(x)-\varphi(y))}{|x-y|^{N+\wp s}}dxdy\\
    &\quad +\vartheta\int\limits_{\mathbb R^N}(|u|^{\frac{N}{s}-2}u+|u|^{q-2}u)\varphi
    dx-\int\limits_{\mathbb R^N}\int\limits_{\mathbb
        R^N}\dfrac{F(u(y))f(u(x))\varphi(x)}{|x-y|^{\mu}} dxdy\
        \hbox{ for every }
        \varphi\in \mathcal {W}.
\end{align*}
 Here  $\langle \cdot,\cdot\rangle$ is the
dual pairing between $\mathcal {W}$ and its dual space $\mathcal {W}'$. Consequently, the
(weak) solutions of problem \eqref{Jeta} are also the critical points of $\mathcal J_{\vartheta}$ in~$\mathcal {W}$.
\begin{lemma}\label{lm2}
Suppose that conditions $(\mathcal {Z}_1)-(\mathcal {Z}_2)$, $(f_1),$ and $(f_5)$ hold. Then there exist constants  $t_0, \rho_0>0$ such that
    $\mathcal E_{\vartheta}(u)\ge \rho_0,$ for every $u\in \mathcal {W},$ with $\|u\|_{\vartheta,\mathcal {W}}=t_0.$
\end{lemma}
\begin{proof}
By condition $(f_1)$, with
$$q_1\ge q> \frac{N}{s}> \frac{2N-\mu}{2s},\quad
q_2\ge\frac{N}{s},$$
there exist $a_1>0, a_2>0$ such that
$$ f'(t)\le a_1|t|^{q_1-2}+a_2\mathcal {H}_{N,s}(\beta_0|t|^{N/(N-s)})|t|^{q_2-2}\
\hbox{ for every } t\in\mathbb R.
$$
 This implies that
$$ f(t)\le a_1|t|^{q_1-1}+a_2\mathcal {H}_{N,s}(\beta_0|t|^{N/(N-s)})|t|^{q_2-1}\
\hbox{ for every }
 t\in\mathbb R.
$$
    Therefore, we get
    \begin{align}\label{F}
        |F(t)|\le a_1|t|^{q_1}+a_2|t|^{q_2}\mathcal H_{N,s}(\beta_0|t|^{N/(N-s)})\
       \hbox{ for every }
       t\in\mathbb R.
    \end{align}
    By Lemma \ref{lm2n}, we obtain that
    \begin{align}\label{d1}
 \int\limits_{\mathbb R^N}\int\limits_{\mathbb R^N}\dfrac{
 |F(u(y))F(u(x))|}{|x-y|^{\mu}}dxdy\le C(r,N,\mu)||F(u)||_{L^{\frac{2N}{2N-\mu}}(\mathbb R^N)}^{2}\
 \hbox{ for every }
 u\in \mathcal {W}.
 \end{align}
Moreover, \eqref{F} yields that
\begin{align}\label{d2a}
 \|F(u)\|_{L^{\frac{2N}{2N-\mu}}(\mathbb R^N)}  \le
 a_1\|u^{q_1}\|_{L^{\frac{2N}{2N-\mu}}(\mathbb R^N)}
 +a_2\||u|^{q_2}\mathcal {H}_{N,s}(\beta_0|u|^{N/(N-s)})\|_{L^{\frac{2N}{2N-\mu}}(\mathbb
 R^N)}.
 \end{align}
Applying the H\"older inequality with $\kappa>1$ and $\kappa'>1$
close to $1$, ${1}/{\kappa}+{1}/{\kappa'}=1$, and using the
arguments from the proof of Li and Yang~\cite[Lemma~2.3 ]{LY},
we can show that
for every
$$\mathfrak l>\frac{2N \kappa'}{2N-\mu},$$
there exists a constant $C(\mathfrak l)>0$ such that
\begin{align}\label{b8a1}
\Big(\mathcal H_{N,s}(\beta_0|t|^{N/(N-s)})\Big)^{\frac{2N\kappa'}{2N-\mu}}\le
C(\mathfrak l)\mathcal H_{N,s}(\mathfrak l \beta_0|t|^{N/(N-s)})\
\hbox{ for every }
t\in\mathbb R.
\end{align}
 Hence, \eqref{b8a1} implies that
for every $u\in \mathcal {W},$
    \begin{align}\label{d3}
&\||u|^{q_2}\mathcal {H}_{N,s}(\beta_0|u|^{N/(N-s)})\|_{{\frac{2N}{2N-\mu}}(\mathbb R^N)}
=        \left(\int\limits_{\mathbb R^N}(|u|^{q_2}\mathcal H_{N,s}(\beta_0|u|^{N/(N-s)}))^{\frac{2N}{2N-\mu}}dx
\right)^{\frac{2N-\mu}{2N}} \notag \\
&\qquad\le        \|u\|_{L^{\frac{2N
\kappa q_2}{2N-\mu}}}^{q_2}\left(\int\limits_{\mathbb
R^N}C(\mathfrak l)\mathcal {H}_{N,s}( \mathfrak l
\beta_0|u|^{N/(N-s)})dx\right)^{\frac{2N-\mu}{2N}}.
\end{align}
    We apply Lemma  \ref{lm1}, taking $\|u\|_{\vartheta,\mathcal {W}}$  small enough, and get
    \begin{align}\label{cta}
        \mathfrak l \beta_0\|u\|_{\vartheta,W^{s,p}(\mathbb R^N)}^{N/(N-s)}\le \mathfrak l \beta_0\|u\|_{\vartheta,\mathcal {W}}^{N/(N-s)}<\alpha_{*},
    \end{align}
    hence
 \begin{equation} \label{d4}
        \int\limits_{\mathbb R^N}
                \mathcal {H}_{N,s}(\mathfrak l \beta_0|u|^{N/(N-s)})dx \\
         =\int\limits_{\mathbb R^N}\mathcal {H}_{N,s}
        \Big(\mathfrak l \beta_0\|u\|_{\vartheta,W^{s,p}(\mathbb{R}^N)}^{N/(N-s)}
        \left(\frac{|u|}{\|u\|_{\vartheta,W^{s,p}(\mathbb{R}^N)}}\right)^{N/(N-s)}\Big)dx<\infty.
    \end{equation}
Together with \eqref{d1}-\eqref{d4}, assuming $\|u\|_{\vartheta,\mathcal {W}}$ to be small enough, we conclude that there exist appropriate constants $\mathfrak h_1>0$ and $\mathfrak h_2>0$ such that
    \begin{align}\label{d5a}
        \int\limits_{\mathbb R^N}\int\limits_{\mathbb R^N}\frac{
        |F(u(y))F(u(x))|}{|x-y|^{\mu}}dxdy\le
        \mathfrak h_1\|u\|_{L^{\frac{2Nq_1}{2N-\mu}}(\mathbb
            R^N)}^{2q_1}+\mathfrak h_2\|u\|_{L^{\frac{2N\kappa q_2}{2N-\mu}}}^{2q_2}.
    \end{align}
Thus, by \eqref{Jeta} and the continuity of the embeddings $\mathcal {W}\hookrightarrow W^{s,N/s}(\mathbb R^N)
\hookrightarrow L^{t}(\mathbb R^N),$
for every $t\ge {N}/{s},$
we obtain that

\begin{align}\label{d5}
\mathcal {E}_{\vartheta}(u)
&\nonumber \ge
\frac{s}{N}\|u\|_{\vartheta,W^{s,p}(\mathbb
R^N)}^{N/s}+\frac{1}{q}\|u\|_{\vartheta,W^{s,q}(\mathbb
R^N)}^{q}
\nonumber
-\mathfrak h_1A_{\frac{2Nq_1}{2N-\mu},\vartheta}^{-2q_1}\|u\|_{\vartheta,\mathcal {W}}^{2q_1}
-\mathfrak h_2A_{\frac{2N\kappa q_2}{2N-\mu},\vartheta}^{-2q_2}\|u\|_{\vartheta,\mathcal {W}}^{2q_2}\\
&\ge \frac{2^{1-q}}{q}||u||_{\vartheta,\mathcal {W}}^{q}-\mathfrak h_1A_{\frac{2Nq_1}{2N-\mu},
\vartheta}^{-2q_1}\|u\|_{\vartheta,\mathcal {W}}^{2q_1}
-\mathfrak h_2A_{\frac{2N\kappa q_2}{2N-\mu},\vartheta}^{-2q_2}\|u\|_{\vartheta,\mathcal {W}}^{2q_2}
\end{align}
for $\|u\|_{\vartheta,\mathcal {W}}$ small enough.
Let
$$\mathscr{C}(t)=\frac{2^{1-q}}{q}-\mathfrak a_1A_{\frac{2Nq_1}{2N-\mu},\vartheta}^{-2q_1}
t^{2q_1-\frac{N}{s}}-\mathfrak h_2A_{\frac{2N\kappa
q_2}{2N-\mu},\vartheta}^{-2q_2}t^{2q_2-\frac{N}{s}},
\quad t\ge 0.$$
We claim that there exists $t_0>0$ so small that
$$\mathscr{C}(t_0)\ge \frac{2^{1-q}}{2q}=\mathscr{C}_0.$$
Clearly, $\mathscr{C}$ is continuous
in $\mathbb R^+_0$ and
$\lim\limits_{t\to 0^{+}}\mathscr{C}(t)=2\mathscr{C}_0,$ so there exists $t_0$ such that
    $\mathscr{C}(t)\ge \mathscr{C}_0,$ for every $t\in [0,t_0].$
We take $t_0$ even smaller, if necessary,  so that
$\|u\|_{\vartheta,\mathcal {W}}=t_0$ satisfies \eqref{cta}. This proves the claim.
Hence $ \mathcal {E}_{\vartheta}(u)\ge \mathscr{C}_0 t_0^{q}=\rho_0,$ for every $u\in \mathcal {W}$, with
$\|u\|_{\vartheta,\mathcal {W}}=t_0$. This completes the proof of Lemma \ref{lm2}.
\end{proof}
In the sequel, $\mathscr{A}_\mu: \mathcal {W}\to\mathbb R$ will denote the
functional
\begin{align}\label{defKmu}
\mathscr{A}_\mu(u)=\frac{1}{2}\int\limits_{\mathbb
R^N}\mathcal {L}_\mu(u)(x)F(u(x))dx,
\end{align}
where $\mathcal {L}_\mu$ is given in~\eqref{Jeta}.

\begin{lemma}\label{lm3}
Suppose that conditions  $(\mathcal {Z}_1), (\mathcal {Z}_2),$ and  $(f_4)$ hold. Then there exists a nonnegative
function $v\in  C_0^{\infty}(\mathbb R^N)$, with
$\|v\|_{\vartheta,\mathcal {W}}>t_0,$ such that $\mathcal {E}_{\vartheta}(v)<0,$ where $t_0>0$ is
the number given by Lemma~$\ref{lm2}$.
\end{lemma}

\begin{proof}
Fix $u_0\in C_0^{\infty}(\mathbb R^N)\setminus \{0\}$, with
$u_0\ge 0$ in $\mathbb R^N$.
Set
$$\mathcal {H}_{\mu}(t)=\mathscr{A}_\mu(tu_0/\|u_0\|_{\vartheta,\mathcal {W}})  \ \ \mbox{for every}
 \ \  t>0,$$
where $\mathscr{A}_\mu$ is defined in \eqref{defKmu}. Condition
$(f_4)$ gives
\begin{align*}
\mathcal {H}_\mu'(t)&=\mathcal {H}_\mu'({tu_0}/{\|u\|_{\vartheta,\mathcal {W}}})\frac{u_0}{\|u_0\|_{\vartheta,\mathcal {W}}}\\
&=\int\limits_{\mathbb
R^N}\Big[{|x|^{-\mu}}\ast
F\left({tu_0}/{\|u_0\|_{\vartheta,\mathcal {W}}}\right)\Big]
f\left({tu_0}/{\|u_0\|_{\vartheta,\mathcal {W}}}\right)\dfrac{u_0}{\|u_0\|_{\vartheta,\mathcal {W}}}dx
>\frac{\theta}{t}\mathcal {H}_\mu(t).
\end{align*}
Integrating the above inequality on $[1,t\|u_0\|_{\vartheta,\mathcal {W}}]$,
with $t>{1}/{\|u_0\|_{\vartheta,\mathcal {W}}},$ we get
$$\mathcal {H}_\mu(t\|u_0\|_{\vartheta,\mathcal {W}})\ge
\mathcal {H}_\mu(1)(t\|u_0\|_{\vartheta,\mathcal {W}})^{\theta}$$ which implies  that
$$\mathscr{A}_\mu(tu_0)\ge \mathscr{A}_\mu\Big(\frac{u_0}{\|u_0\|_{\vartheta,\mathcal {W}}}\Big)\|u_0\|_{\vartheta,\mathcal {W}}^{\theta}t^{\theta}.$$
Therefore, we have
\begin{align*}
        \mathcal {E}_{\vartheta}(tu_0)&=\frac{st^{N/s}}{N}\|u_0\|_{\vartheta,W^{s,p}(\mathbb R^N)}^{N/s}+\frac{t^q}{q}\|u_0\|_{\vartheta,W^{s,q}(\mathbb R^N)}^{q}-\int\limits_{\mathbb R^N}\mathcal {L}_\mu (tu_0)(x)F(tu_0)dx\\
        &\le\frac{st^{N/s}}{N}\|u_0\|_{\vartheta,W^{s,p}(\mathbb R^N)}^{N/s}+\frac{t^q}{q}\|u_0\|_{\vartheta,W^{s,q}(\mathbb R^N)}^{q}-\mathcal {L}_\mu \Big(\frac{u_0}{\|u_0\|_{\vartheta,\mathcal {W}}}\Big)\|u_0\|_{\vartheta}^{\theta}t^{\theta}
\end{align*}
    for every $t>{1}/{\|u_0\|_{\vartheta,\mathcal {W}}},$  choosing $v=tu_0$ and $t$ large
    enough. This completes the proof of Lemma \ref{lm3}.
\end{proof}

Lemmas \ref{lm2} and \ref{lm3} show that $\mathcal {E}_{\vartheta}$
satisfies the geometric conditions
of the  Mountain Pass
Theorem, therefore there exists a Palais-Smale sequence
$\{\mathfrak u_n\}_n\subset \mathcal {W}$ for $\mathcal {E}_{\vartheta}$ at level $c_{\vartheta}$,
briefly $(PS)_{c_\vartheta}$, that is,
\begin{equation}\label{n1}\begin{gathered}\mathcal {E}_{\vartheta}(\mathfrak u_n)\to c_{\vartheta}\quad\text{and}\quad \mathcal {E}'_{\vartheta}(\mathfrak u_n)\to 0
\mbox{ in $\mathcal {W}',$ as }n\to\infty,\\
c_{\vartheta}=\inf_{\zeta\in \Pi}\max_{t\in [0,1]}\mathcal {E}_{\vartheta}(\zeta(t)),
\end{gathered}\end{equation}
where $\Pi=\{\zeta\in C([0,1], \mathcal {W})\,:\, \zeta(0)=0,\;\;
\mathcal {E}_{\vartheta}(\zeta(1))<0\}.$
We denote Nehari manifold $\mathscr{M}_{\vartheta}$ related to $\mathcal {E}_{\vartheta}$ by
\[\mathscr{M}_{\vartheta} = \left\{ {u \in \mathcal {W}\backslash\left\{ 0 \right\}:
\langle {{\mathcal {E}_{\vartheta}'}\left( u \right), u}\rangle  = 0} \right\}.\]

Let us define
\begin{equation}\label{T1}
\mathcal {T}_\vartheta^+:=\{u\in \mathcal {W}: |\text{supp}(u^{+})|>0\}
\end{equation}
and $\mathbb{S}_{Z_0}^+=\mathbb{S}_{Z_0}\cap \mathcal {T}_\vartheta^+,$ where $\mathbb{S}_{\vartheta}$  is the unit sphere in
$\mathcal {W}.$  We know that $\mathcal {T}_\vartheta^+$  is an
open subset of $\mathcal {W}.$

Invoking  the fact above and the definition of $\mathbb{S}_{\vartheta}^+,$
we obtain that $\mathbb{S}_{\vartheta}^+$ is an incomplete $C^{1,1}$-manifold of codimension 1
modelled on  $\mathcal {W}$ and contained in $\mathcal {T}_\vartheta^+.$
Therefore, $\mathcal {W}= T_u\mathbb{S}_{Z_0}^+\oplus \mathbb{R}u$
for every $u\in\mathbb{S}_{Z_0}^+,$
where
\begin{align*}
T_u\mathbb{S}_{\vartheta}^+ & =\{v\in W:
\sum_{\wp\in \{p,q\}}\iint\limits_{\mathbb R^{2N}}\frac{|u(x)-u(y)|^{\wp-2}(u(x)-u(y))(\varphi(x)-\varphi(y))}{|x-y|^{N+\wp s}}dxdy \\
& \quad \quad +\int\limits_{\mathbb R^N}\vartheta(|u|^{p-2}u+|u|^{q-2}u)\varphi dx=0\}.
\end{align*}
Since $f$ is only  continuous, the following result plays an important role in overcoming the nondifferentiability of  $\mathscr{M}_{\vartheta}$ and the
incompleteness of $\mathbb{S}_{\vartheta}.$
\begin{lemma}\label{lm7}
Suppose that conditions $(\mathcal {Z}_1)-(\mathcal {Z}_2)$ and $(f_1)-(f_5)$ hold. Then
\begin{itemize}
\item[$(i)$] For every $u\in \mathcal {T}_\vartheta^+$
and
 $\mathcal {B}_u: [0,\infty)\rightarrow \mathbb{R},$
 defined
as $\mathcal {B}_u(t):= \mathcal {E}_\vartheta(tu),$
 there exists a unique $t_u>0$
such that $\mathcal {B}'_u(t)>0$ in $(0,t_u)$ and $\mathcal {B}'_u(t)<0$ on $(t_u,+\infty).$

\item[$(ii)$]  $\mathscr{M}_\vartheta$ is bounded away from 0 and $\mathscr{M}_\vartheta$ is closed in $\mathcal {W}.$ There exists $\tau>0$ independent on $u,$ such that $t_u\geq \tau,$ for every $u\in \mathbb{S}_\vartheta.$ Moreover, for
each compact set $\mathcal {K}\subset \mathbb{S}_\vartheta,$ there exists $C_{\mathcal {K}}$ such that $t_u\leq C_\mathcal {K},$ for every $u\in \mathcal {K}.$

\item[$(iii)$] The map $\widetilde{m}_{\vartheta}: \mathcal {T}_\vartheta^+\rightarrow\mathscr{M}_{\vartheta},$
given by $\widetilde{m}_{\vartheta}(u)=t_uu,$ is continuous and $m_{\vartheta}:=\widetilde{m}_{\vartheta}|_{\mathbb{S}_{\vartheta}^+}$ is a homeomorphism between $\mathbb{S}_{\vartheta}^+$ and $\mathscr{M}_{\vartheta},$
and  $m_{\vartheta}^{-1}(u)=\frac{u}{\|u\|_{\vartheta,\mathcal {W}}}.$
\end{itemize}
\end{lemma}

\begin{proof}
For every $u\in \mathcal {T}_\vartheta^+$ and $t\geq0,$
$$\mathcal {B}_u(t)=\mathcal {E}_\vartheta(tu)=\frac{t^p}{p}\|u\|_{\vartheta,W^{s,p}(\mathbb R^N)}^{p}+\frac{t^q}{q}\|u\|_{\vartheta,W^{s,q}(\mathbb R^N)}^{q}
-\int\limits_{\mathbb R^N}\int\limits_{\mathbb R^N}\frac{
        |F(tu(y))F(tu(x))|}{|x-y|^{\mu}}dxdy.$$
By \eqref{d5a}, \eqref{d5}  and the continuity of the embeddings $\mathcal {W}\hookrightarrow W^{s,N/s}(\mathbb R^N)\hookrightarrow L^{t}(\mathbb R^N)$  for every $t\in[\frac{N}{s},+\infty),$ we get
\begin{align}\label{d6}
\mathcal {B}_u(t) & \nonumber =\mathcal {E}_\vartheta(tu)
 \nonumber\geq\frac{t^p}{p}\|u\|_{\vartheta,W^{s,p}(\mathbb R^N)}^{p}
+\frac{t^q}{q}\|u\|_{\vartheta,W^{s,q}(\mathbb R^N)}^{q}
- \mathfrak a_1\|u\|_{L^{\frac{2Nq_1}{2N-\mu}}(\mathbb
            R^N)}^{2q_1}-D\|u\|_{L^{\frac{2N\kappa q_2}{2N-\mu}}}^{2q_2} \\
& \geq \frac{t^p}{p}\|u\|_{\vartheta,W^{s,p}(\mathbb R^N)}^{p}
+\frac{t^q}{q}\|u\|_{\vartheta,W^{s,q}(\mathbb R^N)}^{q}
-t^{2q_1}\mathfrak
a_1A_{\frac{2Nq_1}{2N-\mu},\vartheta}^{-2q_1}\|u\|_{\vartheta,\mathcal
{W}}^{2q_1} -t^{2q_2}DA_{\frac{2N\kappa
q_2}{2N-\mu},\vartheta}^{-2q_2}\|u\|_{\vartheta,\mathcal {W}}^{2q_2}
\end{align}
which yields $\mathcal {B}_u(t)\rightarrow 0^+,$ as $t\rightarrow0^+.$ Moreover, we obtain that
$\mathcal {B}_u(t)\rightarrow -\infty,$ as $t\rightarrow\infty.$ Therefore, there exists $t_u\in (0,\infty)$
such that $\mathcal {B}_u(t_u)=\max_{t\geq0}\mathcal {B}_u(t).$ Furthermore, $\mathcal {B}'_u(t_u)=0.$
Now we shall verify that $t_u$ is a unique critical point of $\mathcal {B}_u$
in $(0,\infty).$ Arguing by contradiction, suppose that there exist $0<t_1<t_2<\infty$ such that $\mathcal {B}'_u(t_1)=\mathcal {B}'_u(t_2)=0.$ Consequently, we have
$$t_1^p\|u\|_{\vartheta,W^{s,p}(\mathbb R^N)}^{p}+t_1^q\|u\|_{\vartheta,W^{s,q}(\mathbb R^N)}^{q}
=\int\limits_{\mathbb R^N}\int\limits_{\mathbb R^N}\frac{|F(t_1u(y))f(t_1u(x))t_1u(x)|}{|x-y|^{\mu}}dxdy$$
and
$$t_2^p\|u\|_{\vartheta,W^{s,p}(\mathbb R^N)}^{p}+t_2^q\|u\|_{\vartheta,W^{s,q}(\mathbb R^N)}^{q}
=\int\limits_{\mathbb R^N}\int\limits_{\mathbb R^N}\frac{|F(t_2u(y))f(t_2u(x))t_2u(x)|}{|x-y|^{\mu}}dxdy.$$
By two equalities above and $(f_4),$ we obtain
\begin{align*}
 0&<(\frac{1}{t_1^{q-p}}-\frac{1}{t_2^{q-p}})\|u\|_{\vartheta,W^{s,p}(\mathbb R^N)}^{p} \\
& = \int\limits_{\mathbb R^N}\int\limits_{\mathbb R^N}
\frac{1}{|x-y|^{\mu}}\Big|\frac{F(t_1u(y))}{(t_1u(y))^{\frac{q}{2}}}
\frac{f(t_1u(x))}{(t_1u(x))^{\frac{q}{2}-1}}(u(y))^{\frac{q}{2}}(u(x))^{\frac{q}{2}}\Big| dydx  \\
& \quad -\int\limits_{\mathbb R^N}\int\limits_{\mathbb R^N}
\frac{1}{|x-y|^{\mu}}\Big|\frac{F(t_2u(y))}{(t_2u(y))^{\frac{q}{2}}}
\frac{f(t_2u(x))}{(t_2u(x))^{\frac{q}{2}-1}}(u(y))^{\frac{q}{2}}(u(x))^{\frac{q}{2}}\Big| dydx \\
& <0
\end{align*}
which is impossible. Therefore, we have completed the proof of $(i).$

$(ii)$ Let $u\in \mathscr{M}_\vartheta.$ We shall prove  that  the
first part of conclusion is true in the following two cases.

\medskip

$\textbf{Case 1.}$ $\|\mathfrak u_n\|_{\vartheta,W^{s,p}(\mathbb{R}^N)}^{N/(N-s)}>\frac{\alpha_{*}}{ \mathfrak l \beta_0}\sigma^{s/(N-s)}.$

 In this case we are done.

\medskip

$\textbf{Case 2.}$ $\|\mathfrak u_n\|_{\vartheta,W^{s,p}(\mathbb{R}^N)}^{N/(N-s)}<\frac{\alpha_{*}}{ \mathfrak l \beta_0}\sigma^{s/(N-s)}.$

Applying the Trudinger-Moser inequality,
we obtain
\begin{align}\label{nm3o}
\sup_{u\in W^{s,p}(\mathbb R^N)}\int\limits_{\mathbb R^N}\Phi_{N,s}(\alpha_0 |u|^{N/(N-s)})dx<+\infty\
\hbox{ for every }
0\le \alpha<\alpha_{*}.
\end{align}

Using the Hardy-Littlewood-Sobolev inequality again and $(f_3)$, it
follows that
 $$\int\limits_{\mathbb R^N}\int\limits_{\mathbb R^N}\frac{|F(u(y))f(u(x))u(x)|}{|x-y|^{\mu}}dydx
 \leq C\|F(u)\|_{L^{\frac{2N}{2N-\mu}}}\|f(u)u\|_{L^{\frac{2N}{2N-\mu}}}
 \leq C \|f(u)u\|^2_{L^{\frac{2N}{2N-\mu}}}. $$
By $(f_1)$ and $(f_2),$ for every $\varepsilon_{*}>0$ and $q_1\ge q> \frac{N}{s}> \frac{2N-\mu}{2s}$ and
$q_2\ge\frac{N}{s},$
 there exists $C_{q,\varepsilon_{*}}>0$ such that
\begin{align}\label{nm3}
\|f(u)u\|_{L^{\frac{2N}{2N-\mu}}}\le \varepsilon_{*} \||u|^{q_1}\|_{L^{\frac{2N}{2N-\mu}}}+C_{q,\varepsilon_{*}}\||u|^{q_2}\mathcal {H}_{N,s}
(\alpha_0|t|^{N/(N-s)})\|_{L^{\frac{2N}{2N-\mu}}}
\end{align}
for every $t\ge 0.$
Using inequality (\ref{nm3}) and the definition of $A_{\nu,\vartheta}$, there exists a constant $C(\varepsilon_{*})$ such that
\begin{align}\label{ct00.1}
\|f(u)u\|^2_{L^{\frac{2N}{2N-\mu}}}
\le \varepsilon_{*}A_{\frac{2Nq_1}{2N-\mu},
\vartheta}^{-2q_1}\|u\|_{\vartheta,\mathcal {W}}^{2q_1}
+C(q,\varepsilon_{*})\|u\|_{\vartheta,\mathcal {W}}^{2q_2}\
\hbox{ for some }
q_1\ge q, q_2\ge \frac{N}{s}.
\end{align}
 In view of $\langle \mathcal {E}'_{\vartheta}(u),u\rangle=0,$ we get
\begin{align*}
\|u\|_{\vartheta,W^{s,p}(\mathbb R^N)}^{p}+\|u\|_{\vartheta,W^{s,q}(\mathbb R^N)}^{q}
& =\int\limits_{\mathbb R^N}\int\limits_{\mathbb R^N}\frac{|F(u(y))f(u(x))u(x)|}{|x-y|^{\mu}}dydx \\
& \le  \varepsilon_{*}A_{\frac{2Nq_1}{2N-\mu},
\vartheta}^{-2q_1}\|u\|_{\vartheta,\mathcal {W}}^{2q_1}
+C(q,\varepsilon_{*})\|u\|_{\vartheta,\mathcal {W}}^{2q_2}.
\end{align*} Therefore
\begin{equation}\label{e2.14}
\|u\|_{\vartheta, \mathcal {W}}\geq \alpha\ \ \hbox{for every } u\in \mathscr{M}_{\vartheta}.
\end{equation}
For any sequence $\{\mathfrak u_n\}_n\subset\mathscr{M}_{\vartheta},$ such that $\mathfrak u_n\to u$ in $\mathcal {W},$ we have to prove that $u\in\mathscr{M}_{\vartheta}.$ Indeed,
by the fact that $\mathfrak u_n\to u$ in $\mathcal {W},$ we obtain $\|\mathfrak u_n-u\|_{\vartheta,\mathcal {W}}\to 0,$ as $n\to\infty.$ Using the discussion as in
Willem \cite[Lemma A.1]{NW96},
 there exists a subsequence $\{\mathfrak u_n\}_n$ of $\{u_n\}_n$ satisfying
$\|\mathfrak u_{i+1}-\mathfrak u_i\|_{\vartheta,\mathcal {W}}\le 2^{-i}$ for every $i\ge 1.$ We denote
$\mathcal {U}(x):=|\mathfrak u_1(x)|+\sum_{i=1}^{\infty}|\mathfrak u_{i+1}(x)-\mathfrak u_i(x)|.$
Together with the fact that $\mathfrak u_n\to u$ in $\mathcal {W},$ we obtain
$|\mathfrak u_n(x)|\le \mathcal {U}(x)$ for every $x\in \mathbb R^N$ and $|u(x)|\le \mathcal {U}(x)$ in $\mathbb R^N.$ Clearly,
$\|\mathcal {U}\|_{\vartheta,\mathcal {W}}\le \|\mathfrak u_1\|_{\vartheta,\mathcal {W}}+\sum_{i=1}^{\infty}2^{-i}<+\infty,$ hence $\mathcal {U}\in \mathcal {W}.$

Since $\{\mathfrak u_n\}_n$ is a subsequence of $\{\mathfrak u_n\}_n,$  we have
\begin{align}\label{ct1b}
    \|u\|_{\vartheta,W^{s,p}(\mathbb R^N)}^{p}
    + \|u\|_{\vartheta,W^{s,q}(\mathbb R^N)}^{q}
    =\lim_{n\to\infty}\int\limits_{\mathbb R^N}\int\limits_{\mathbb R^N}\frac{|F(\mathfrak u_n(y))f(\mathfrak u_n(x))\mathfrak u_n(x)|}{|x-y|^{\mu}}dydx.
\end{align}
We shall now prove
\begin{align}\label{ct1}
    \lim_{n\to\infty}\int\limits_{\mathbb R^N}\int\limits_{\mathbb R^N}\frac{|F(\mathfrak u_n(y))f(\mathfrak u_n(x))\mathfrak u_n(x)|}{|x-y|^{\mu}}dydx
    =\int\limits_{\mathbb R^N}\int\limits_{\mathbb R^N}\frac{|F(\mathfrak u(y))f(\mathfrak u(x))\mathfrak u(x)|}{|x-y|^{\mu}}dydx.
\end{align}
We have
\begin{align*}
    &\Big|\frac{|F(\mathfrak u_n(y))f(\mathfrak u_n(x))\mathfrak u_n(x)|}{|x-y|^{\mu}}-\frac{|F(\mathfrak u(y))f(\mathfrak u(x))\mathfrak u(x)|}{|x-y|^{\mu}}\Big|
    \le 2\frac{|F(\mathcal {U}(y))f(\mathcal {U}(x))\mathfrak u(x)|}{|x-y|^{\mu}}.
\end{align*}
Now we shall show that
\begin{align}\label{ctmg0}
    \frac{|F(\mathcal {U}(y))f(\mathcal {U}(x))\mathcal {U}(x)|}{|x-y|^{\mu}}\in L^1(\mathbb R^N).
\end{align}
By
Zhang et al. \cite[Lemma 2.4]{ZXN23}, we have
\begin{align}\label{nm3oa}
    \int\limits_{\mathbb R^N}\mathcal {H}_{N,s}(\alpha_0 |\mathcal {U}|^{N/(N-s)})dx<+\infty.
\end{align}
Together with \eqref{nm3}-\eqref{ct00.1}, we deduce
(\ref{ctmg0}). Therefore,
$$\int\limits_{\mathbb R^N}\int\limits_{\mathbb R^N}\frac{|F(\mathfrak u_n(y))f(\mathfrak u_n(x))\mathfrak u_n(x)|}{|x-y|^{\mu}}dydx-\int\limits_{\mathbb R^N}\int\limits_{\mathbb R^N}\frac{|F(u(y))f(u(x))u(x)|}{|x-y|^{\mu}}dydx\in L^1(\mathbb{R}^N).$$
Furthermore,
$$\int\limits_{\mathbb R^N}\int\limits_{\mathbb R^N}\frac{|F(\mathfrak u_n(y))f(\mathfrak u_n(x))\mathfrak u_n(x)|}{|x-y|^{\mu}}dydx-\int\limits_{\mathbb R^N}\int\limits_{\mathbb R^N}\frac{|F(u(y))f(u(x))u(x)|}{|x-y|^{\mu}}dydx\rightarrow 0, $$ pointwisely on $\mathbb R^N$ outside a set of measure zero. By the Dominated Convergence Theorem, we obtain that \eqref{ct1} is true. By $\mathfrak u_n\to u$ in $\mathcal {W}(\mathbb R^N),$ we have
$$\|u\|_{\vartheta,W^{s,p}(\mathbb R^N)}^{p}
+\|u\|_{\vartheta,W^{s,q}(\mathbb R^N)}^{q}=\int\limits_{\mathbb R^N}\int\limits_{\mathbb R^N}\frac{|F(u(y))f(u(x))u(x)|}{|x-y|^{\mu}}dydx $$
which yields $u\in \mathscr{M}_{\vartheta}.$

In the sequel, we shall verify that  the second part of the conclusion is also true.

Applying $(i),$  there exists $t_u>0$ such that
$t_uu\in\mathscr{M}_{\vartheta}$ for every $u\in \mathbb{S}_{\vartheta}.$
Therefore, it follows from \eqref{e2.14} that  $t_u\geq\alpha.$
We shall argue by contradiction that
$\mathfrak u_n\in \mathcal {K}$ satisfies $t_n:= t_{\mathfrak u_n}\rightarrow\infty.$ Due to the compactness of  $\mathcal {K},$ we may suppose that
$\mathfrak u_n\rightarrow u$ in $\mathcal {W}.$ Then $u\in\mathcal {K}\subset\mathbb{S}_{\vartheta}.$ By $(f_4)$, we obtain
\begin{align*}
{\mathcal {E}_{\vartheta}(t_n\mathfrak u_n)}& = \frac{1}{p}t_n^p\left\| {{\mathfrak u_n}} \right\|_{\vartheta,W^{s,p}(\mathbb R^N)}^{p}
+\frac{1}{q}t_n^p\left\| {{\mathfrak u_n}} \right\|_{\vartheta,W^{s,p}(\mathbb R^N)}^{p}
 - \frac{1}{2}\int_{{\mathbb{R}^N}} {\int_{{\mathbb{R}^N}} {\frac{{F({t_n}{\mathfrak u_n}(y))}}{{|x - y{|^\mu }}}} } F({t_n}{\mathfrak u_n}(x))dxdy \\
&\geq\frac{1}{p}t_n^p\left\| {{\mathfrak u_n}} \right\|_{\vartheta,W^{s,p}(\mathbb R^N)}^{p}
+\frac{1}{q}t_n^p\left\| {{\mathfrak u_n}} \right\|_{\vartheta,W^{s,q}(\mathbb R^N)}^{q}
-\gamma_1^2t_n^{2\theta}\||\mathfrak u_n|^\theta\|_{L^{\frac{2N}{2N-\mu}}(\mathbb{R}^N)}^2\rightarrow-\infty
\end{align*} due to $\theta>q$. However, since $t_n\mathfrak u_n\in \mathscr{M}_{\vartheta},$ we have
\begin{align*}
\mathcal {E}_{\vartheta}|_{\mathscr{M}_{\vartheta}}(t_n\mathfrak u_n)& =\frac{1}{p}t_n^p\left\| {{\mathfrak u_n}} \right\|_{\vartheta,W^{s,p}(\mathbb R^N)}^{p}
+\frac{1}{q}t_n^p\left\| {{\mathfrak u_n}} \right\|_{\vartheta,W^{s,q}(\mathbb R^N)}^{q}
-\frac{1}{2}\int_{\mathbb{R}^N}\int_{\mathbb{R}^N}\frac{F(\mathfrak u_n(y))F(\mathfrak u_n(x))}{|x-y|^\mu}dydx \\
& \leq\int_{\mathbb{R}^N}\int_{\mathbb{R}^N}\frac{F(t_n\mathfrak u_n(y))}{|x-y|^\mu}
\Big[\frac{1}{p}f(t_n\mathfrak u_n(x))t_n\mathfrak u_n(x)-\frac{1}{2}F(t_n\mathfrak u_n(x))\Big]dydx\geq 0
\end{align*}
which is impossible.

\item[$(iii)$] By $(i)-(ii)$ and the arguments from the proof of
Szulkin and Weth \cite[Proposition 3.1]{LQ}, we obtain $(iii).$ This completes the proof of Lemma \ref{lm7}.
\end{proof}
\begin{remark}\label{rem1}
By Lemma \ref{lm7}, the least energy $c_{Z_0}$ satisfies the following equality:
\begin{equation}\label{e2.19}
c_{\vartheta}=\inf_{u\in \mathscr{M}_{\vartheta}}\mathcal {E}_{\vartheta}(u)
=\inf_{u\in \mathcal {W}\backslash\{0\}}\max_{t>0}\mathcal {E}_{\vartheta}(tu)
=\inf_{u\in \mathbb{S}_{\vartheta}}\max_{t>0}\mathcal {E}_{\vartheta}(tu).
\end{equation}
Considering the functional $\Phi_{\vartheta}:\mathcal {S}_{\vartheta}\rightarrow\mathbb{R}$
given by
\begin{equation}\label{3b}
\Phi_{\vartheta}(\omega):= \mathcal
{E}_{\vartheta}(m_{\vartheta}(\omega))
\end{equation}
 similarly to
 Szulkin and Weth \cite [Corollary 3.3]{LQ}, we have:
\end{remark}
\begin{lemma}\label{lm7.1}
Suppose that conditions $(\mathcal {Z}_1)-(\mathcal {Z}_2)$ and $(f_1)-(f_5)$ hold. Then the following statements are true:
\begin{itemize}
\item[$(i)$] If $\{\mathfrak u_n\}_n$ is a $(PS)_{c_{\vartheta}}$ sequence for $\Phi_{\vartheta},$ then $\{m_{\vartheta}(\mathfrak u_n)\}$ is
a $(PS)_{c_{\vartheta}}$ sequence for $\mathcal {E}_{\vartheta}.$ If $\{\mathfrak u_n\}_n\subset\mathscr{M}_{\vartheta}$ is a bounded $(PS)_{c_{\vartheta}}$ sequence for $\mathcal {E}_{\vartheta},$
then $\{m_{\vartheta}^{-1}(\mathfrak u_n)\}_n$ is a $(PS)_{c_{\vartheta}}$ sequence for $\Phi_{\vartheta}.$

\item[$(ii)$]  $u$ is a critical point of $\Phi_{\vartheta}$ if and only if $m_{\vartheta}(u)$ is a nontrivial critical point of $\mathcal {E}_{\vartheta}.$ Moreover, $\inf_{\mathscr{M}_{\vartheta}}\mathcal {E}_{\vartheta}=\inf_{\mathcal {S}_{\vartheta}}\Phi_{\vartheta}.$
\end{itemize}
\end{lemma}
\begin{lemma}\label{lma}(see Liang et al. \cite[Lemma 5]{SLN23})
Let $(f_1)$ be satisfied. We suppose that $\{\mathfrak u_n\}_n$
is a sequence in $W^{s,p}(\mathbb{R}^N)$  such that
$$\limsup\limits_{n\to\infty} \|\mathfrak u_n\|_{\vartheta,W^{s,p}(\mathbb{R}^N)}^{N/(N-s)}<\frac{\alpha_{*}}{ \mathfrak l \beta_0}\sigma^{s/(N-s)} 
 \ \ \mbox{for some} \ \  \mathfrak l>1,$$
where $\sigma=\min\{1,\vartheta\}.$  Then
there exists $C_0>0$ such that
$$ \Big|{|x|^{-\mu}}*F(\mathfrak u_n)\Big|\le C_0\
\mbox{ for every } n. $$
\end{lemma}

\begin{lemma}(see Molica Bisci et al. \cite[Lemma 4]{NVT1})\label{lm3b}
Let $\varsigma\in[{N}/{s},\infty).$
    If $\{\mathfrak u_n\}_n$ is a bounded sequence in $W^{s,p}(\mathbb R^N)$ and
    $$ \lim\limits_{n\to\infty}\sup_{y\in\mathbb R^N}\int\limits_{B_R(y)}|\mathfrak u_n(x)|^{\varsigma}dx=0$$
    for some $R>0,$ then $\mathfrak u_n\to 0$ in $L^{\nu}(\mathbb R^N),$ for every $\nu\in
    (\varsigma,\infty).$
\end{lemma}

As in the proof of  Molica Bisci et al.~\cite[Theorem 7]{NVT1}, we
get
\begin{lemma}\label{lm3c}
Suppose that conditions $(\mathcal {Z}_1)-(\mathcal {Z}_2)$ and $(f_1)$ hold. Then
$$\limsup\limits_{t\to0^{+}}f(t)t^{1-\varsigma}=0$$
for some $\varsigma\ge N/s.$
Let $\{\mathfrak u_n\}_n\subset \mathcal {W}$ be weakly convergent
to $0$ and such that
$$\limsup\limits_{n\to\infty}\|\mathfrak u_n\|_{\vartheta,W^{s,p}(\mathbb{R}^N)}^{N/(N-s)}< \frac{\alpha_{*}}{\mathfrak l \beta_0}\sigma^{s/(N-s)},
$$
where  $\beta_{0}$, ${\alpha_*}$ are given in $(f_1)$ and
\eqref{alpha*},
 $\sigma=\min\{1,\vartheta\},$ and
 $\mathfrak l>1$ is a
suitable constant.
If there exists $R>0$ such that
$$\liminf\limits_{n\to\infty}\sup\limits_{y\in\mathbb
R^N}\int\limits_{B_R(y)}|\mathfrak u_n|^{\varsigma}dx=0,$$
then
$$ \int\limits_{\mathbb R^N}[{|x|^{-\mu}}\ast F(\mathfrak u_n)]f(\mathfrak u_n)\mathfrak u_n\to 0\quad \text{and}\quad \int\limits_{\mathbb R^N}[{|x|^{-\mu}}\ast F(\mathfrak u_n)]F(\mathfrak u_n)\to 0\
\hbox{ as }
n\to\infty.$$

\end{lemma}
\begin{lemma}\label{lm3a}(see Liang et al. \cite[Lemma 2.6]{Liang1})
Suppose  that conditions  $(f_3)$ and $(f_4)$ hold. Then there exists a constant $C_{\gamma_1}$ such that $\rho_0\le c_{\vartheta}\le C_{\gamma_1}$,
where $\rho_0$ is the number determined by  Lemma~$\ref{lm2},$
$\gamma_1$ is the constant given in $(f_4)$, $$C_{\gamma_1}=a\left(1-\frac{N}{2\theta s}\right)
    \left(\frac{a N}{2\theta sb}\right)^{N/(2\theta s-N)},$$ and
$$a= \frac{s}{N}\|u\|_{\vartheta,W^{s,p}(\mathbb{R}^N)}^{\frac{N}{s}}
+ \frac{1}{q}\|u\|_{\vartheta,W^{s,q}(\mathbb{R}^N)}^{q}.$$
\end{lemma}
\begin{proposition}\label{pro1}
Suppose that conditions $(\mathcal {Z}_1)-(\mathcal {Z}_2)$ and $(f_1)-(f_5)$ hold. Then
problem \eqref{ct1a} has a nontrivial nonnegative (weak)
solution.
\end{proposition}
\begin{proof}
By Lemmas \ref{lm2} and~\ref{lm3},  there exists a $(PS)_{c_\vartheta}$
sequence $\{\mathfrak u_n\}_n\subset \mathcal {W}$,  satisfying \eqref{n1}.
We shall divide the proof into two steps.

\medskip

$\textbf{Step 1.}$ $\{\mathfrak u_n\}_n$ is a bounded sequence in $\mathcal {W}.$

 Up to a subsequence, we may suppose that $\{\mathfrak u_n\}_n$ is strongly convergent in $\mathcal {W}.$ From \eqref{n1}, we have
\begin{align}\label{n2}
\mathcal {E}_\vartheta(\mathfrak u_n)-\frac{1}{\theta}\langle \mathcal {E}'_{\vartheta}(\mathfrak u_n), \mathfrak u_n\rangle= c_{\vartheta}+o_n(1)+o_n(1)\|\mathfrak u_n\|_{\vartheta, \mathcal {W}}\
\hbox{ as }
n\to\infty,
\end{align}
  where $\theta$ is given in $(f_3).$ Moreover,
\begin{align*}
\mathcal {E}_{\vartheta}(\mathfrak u_n)-\frac{1}{\theta}\langle \mathcal {E}_{\vartheta}'(\mathfrak u_n), \mathfrak u_n\rangle=&\Big(\frac{s}{N}-\frac{1}{\theta}\Big)\|\mathfrak u_n\|_{\vartheta,W^{s,p}(\mathbb R^N)}^{N/s}+\Big(\frac{1}{q}-\frac{1}{\theta}\Big)\|\mathfrak u_n\|_{\vartheta,W^{s,q}(\mathbb R^N)}^{q}\\
&+\int\limits_{\mathbb R^N}\int\limits_{\mathbb
R^N}\frac{F(\mathfrak u_n)}{|x-y|^{\mu}}\Big[\dfrac{1}{\theta}f(\mathfrak u_n)\mathfrak u_n-\frac{1}{2}F(\mathfrak u_n)\Big]
dxdy.
\end{align*}
Condition $(f_2)$ implies that
\begin{align}\label{n4}
\mathcal {E}_{\vartheta}(\mathfrak u_n)-\frac{1}{\theta}&\langle \mathcal {E}'_{\vartheta}(\mathfrak u_n), \mathfrak u_n\rangle\ge
\Big(\frac{s}{N}-\frac{1}{\theta}\Big)\|\mathfrak u_n\|_{\vartheta,W^{s,p}(\mathbb
R^N)}^{N/s}+\Big(\frac{1}{q}-\frac{1}{\theta}\Big)\|\mathfrak u_n\|_{\vartheta,W^{s,q}(\mathbb
R^N)}^{q}.
\end{align}
Combining \eqref{n2} and \eqref{n4}, we get as $n\to\infty,$
\begin{align}\label{n5}
\Big(\frac{s}{N}-\frac{1}{\theta}\Big)\|\mathfrak u_n\|_{\vartheta,W^{s,p}(\mathbb
R^N)}^{N/s}+\Big(\frac{1}{q}-\frac{1}{\theta}\Big)
\|\mathfrak u_n\|_{\vartheta,W^{s,p}(\mathbb R^N)}^{q}\le
c_{\vartheta}+o_n(1)+o_n(1)\|\mathfrak u_n\|_{\vartheta, \mathcal {W}}.
\end{align}
Observing that $\lim\limits_{x\to\infty,\,y\to\infty}\frac{\mathfrak c x^{N/s}+\mathfrak k y^{q}}{x+y}=\infty,$
for fixed numbers $\mathfrak c>0$, $\mathfrak k>0,$
we conclude that $\{\mathfrak u_n\}_n$ is a bounded sequence in $\mathcal {W}.$
Since
$\mathcal {E}_{\vartheta}(\mathfrak u_n)-\frac{1}{\theta}\langle \mathcal {E}_{\vartheta}'(\mathfrak u_n), \mathfrak u_n\rangle\to c_{\vartheta}$ as $n\to \infty,$
it follows that
\begin{align}\label{n6}
\limsup_{n\to\infty}\|\mathfrak u_n\|_{\vartheta, W^{s,p}(\mathbb{R}^N)}^{N/s}\le
\frac{c_{\vartheta}}{\frac{1}{p}-\frac{1}{\theta}}\le
\frac{C_{\gamma_1}}{\frac{1}{p}-\frac{1}{\theta}}
\end{align}
and
\begin{align}\label{n6n}
\limsup\limits_{n\to\infty}\|\mathfrak u_n\|_{\vartheta,
W^{s,q}(\mathbb{R}^N)}^{q}\le
\frac{c_{{\vartheta}}}{\frac{1}{q}-\frac{1}{\theta}}\le
\frac{C_{\gamma_1}}{\frac{1}{q}-\frac{1}{\theta}}.
\end{align}
From \eqref{n6}, we have
\begin{align}\label{cte1new}
\limsup_{n\to\infty} \|\mathfrak u_n\|_{\vartheta, W^{s,p}(\mathbb{R}^N)}^{N/(N-s)}<\frac{\alpha_{*}}{ \mathfrak l \beta_0}\sigma^{s/(N-s)}
\end{align}
for some $\mathfrak l>0$ and $\gamma_1$ large enough such that
\begin{align*}
\frac{C_{\gamma_1}}{\frac{1}{p}-\frac{1}{\theta}}< \Big(\frac{\alpha_*}
{\mathfrak l \beta_0}\Big)^{(N-s)/N}\sigma,
\end{align*}
when $\gamma_1\ge \gamma_{*}$, where
\begin{equation}\label{1a}
\gamma_{*}=\frac{1}{|NB_1(0)|}\sqrt{\frac{2a(N-\mu)}{\mathfrak B(N,N-\mu+1)}}
\end{equation}
 and $\mathfrak B=\int_0^1t^{x-1}(1-t)^{y-1}dt, x>0,y>0$ is  Beta function. This
implies that
\begin{align*}
\gamma_1\ge \frac{1}{|NB_1(0)|}\sqrt{\frac{aN(N-\mu)}{\theta a\mathfrak B(N,N-\mu+1)}}\Bigg[\frac{a(1-
    \frac{N}{2\theta s})}{(\frac{s}{N}-\frac{1}{\theta})\sigma (
    \frac{\alpha_*}{\mathfrak b \beta_0})^{(N-s)/N}}\Bigg]^{\frac{2\theta s-N}{2N}}=\gamma_{**}.
\end{align*}
Hence, \eqref{cte1new} holds for every $\gamma_1\ge \max\{
\gamma_{*},\gamma_{**}\}$.

\medskip

\textbf{Step 2.}~~
We shall show that there exist $R>0,\delta>0$ and a sequence
$\{y_n\}_n\subset\mathbb R^N$ such that
\begin{align}\label{ctabc1}
\liminf_{n\to\infty}\int\limits_{B_R(y_n)}|\mathfrak u_n(x)|^{\wp}dx\ge
\delta>0,\quad\wp\in \{p,q\}.
\end{align}

Arguing by contradiction, we assume that for some $R>0,$ we have
 \begin{align}\label{mn1}
 \lim_{n\to\infty}\sup_{y\in\mathbb R^N}\int\limits_{B_R(y)}|\mathfrak u_n(x)|^{\wp}dx=0,
 \quad {\wp}\in\{p,\,q\}.
 \end{align}
Then Lemma \ref{lm3b} yields that $\mathfrak u_n\to 0$ in $L^{\nu}(\mathbb R^N),$
for every  $\nu>\wp.$ Condition \eqref{mn1}, Lemma~\ref{lm3c}
and the Trudinger-Moser inequality imply that
    \begin{align*}
        \int\limits_{\mathbb R^N}\Big[{|x|^{-\mu}}*F(\mathfrak u_n)\Big]f(\mathfrak u_n)\mathfrak u_ndx\to 0\ \text{as}\; n\to \infty.
    \end{align*}
    Therefore,
    \begin{align*}
        o_n(1)&=\langle \mathcal E_{\vartheta}'(\mathfrak u_n),\mathfrak u_n\rangle=\|\mathfrak u_n\|_{\vartheta,W^{s,p}(\mathbb R^N)}^{p}+\|\mathfrak u_n\|_{\vartheta,W^{s,q}(\mathbb R^N)}^{q}-\int\limits_{\mathbb R^N}\Big[{|x|^{-\mu}}*F(\mathfrak u_n)\Big]f(\mathfrak u_n)\mathfrak u_ndx\\
        &=\|\mathfrak u_n\|_{\vartheta,W^{s,p}(\mathbb R^N)}^{p}+\|\mathfrak u_n\|_{\vartheta,W^{s,q}(\mathbb R^N)}^{q}+o_n(1)\
        \hbox{ as } n\to\infty.
    \end{align*}
     So $\mathfrak u_n\to 0$  in $W^{s,N/s}(\mathbb R^N)\cap W^{s,q}(\mathbb R^N).$
   Passing to the limit as $n\rightarrow\infty,$ we have
    $$\mathcal E_{\vartheta}(\mathfrak u_n)=\frac{\|\mathfrak u_n\|_{\vartheta,W^{s,p}(\mathbb R^N)}^{p}}{p}+\frac{\|\mathfrak u_n\|_{\vartheta,W^{s,q}(\mathbb R^N)}^{q}}{q}-\int\limits_{\mathbb R^N}\Big[{|x|^{-\mu}}*F(\mathfrak u_n)\Big]F(\mathfrak u_n)dx\to 0 $$
which contradicts with the fact that $\mathcal E_{\vartheta}(\mathfrak u_n)\to c_{\vartheta}>0,$
as $n\to\infty$. Thus \eqref{ctabc1} holds.

Put $\mathfrak v_n=\mathfrak u_n(\cdot+y_n),$ then from \eqref{ctabc1} we have
\begin{align}\label{mn2}
\int\limits_{B_R(0)}|\mathfrak v_n|^{\wp}dx\ge \delta/2>0,
\quad \wp\in\{p,q\}\
\hbox{ for some } \delta>0.
\end{align}
  Because $\mathcal E_{\vartheta}$ and $\mathcal E_{\vartheta}'$ are invariant under  translation, we have
$$ \mathcal E_{\vartheta}(\mathfrak v_n)\to c_{\vartheta}\; \ \ \text{and}\; \ \
\mathcal E_{\vartheta}'(\mathfrak v_n)\to 0\;  \ \ \text{in}\; \mathcal {W}'.$$
Since $\|\mathfrak v_n\|_{\vartheta,\mathcal {W}}=\|\mathfrak u_n\|_{\vartheta,\mathcal {W}}$  for every $n$,
then $\{\mathfrak v_n\}_n$ is  also bounded in $\mathcal {W}$ and
\begin{align}\label{cte1newa}
    \limsup_{n\to\infty} \|\mathfrak v_n\|_{\vartheta, W^{s,p}(\mathbb{R}^N)}^{N/(N-s)}=\limsup_{n\to\infty} \|\mathfrak u_n\|_{\vartheta, W^{s,p}(\mathbb{R}^N)}^{N/(N-s)}<\frac{\alpha_{*}}{ \mathfrak l \beta_0}\sigma^{s/(N-s)}.
\end{align}
Hence, choosing a subsequence if necessary, we may assume that there exists $\mathfrak v\in \mathcal {W}$ such
that $\mathfrak v_n\rightharpoonup \mathfrak v$
in $W,$ $\mathfrak v_n\to \mathfrak v$ in $L^{\vartheta}(B_R(0)),$
for every  $\vartheta\in [{N}/{s},\infty)$ and $R>0$,  and $\mathfrak v_n\to \mathfrak v$ a.e. in $\mathbb R^N.$ Clearly, \eqref{mn2} implies that
$$\int\limits_{B_R(0)}|\mathfrak v|^{\wp}dx\ge \delta/2>0,\quad \wp\in\{p,q\},$$
hence, $\mathfrak v\not\equiv 0.$

Arguing  similarly  to the proof of  Thin et al.
\cite[ Lemma~13]{TTL23}, we get that $\mathcal E_{\vartheta}'(\mathfrak v)=0$
and  $\mathfrak v$ is indeed a ground state solution of problem
\eqref{ct1a}.
This completes the proof of Proposition \ref{pro1}.
\end{proof}

\section{The auxiliary problem \eqref{m4}}\label{s4}
\def\theequation{4.\arabic{equation}}
\setcounter{equation}{0}

Using the transformation $x\mapsto \varepsilon x,$ problem
\eqref{pr2} can be  rewritten as follows
\begin{align*}\label{m4}
    (-\Delta)_{p}^{s}u+(-\Delta)_q^su+Z(\varepsilon x)(|u|^{p-2}u+|u|^{q-2}u)=[{|x|^{-\mu}}\ast
    F(u)]f(u). \tag{$\mathcal {Q}_{\varepsilon}$}
\end{align*}
Inspired by the work of  del Pino and Felmer~\cite{Pi1},
we introduce a penalized function which will play  an essential role to obtain our main results.
In general, we assume that
$0\in\Omega$ and $Z(0)=Z_0.$  Fix $\hbar_0>0.$
We define
\begin{align*}
 \hat{f}(t):=
\begin{cases}
&f(t)\;\hspace{0.49cm}  \ \ \text{if}\; t\le a,\\
&\frac{Z_0}{\hbar_0}t^{q-1}\; \ \  \text{if}\; t>a,
\end{cases}
\end{align*}
and
\begin{align*}
 \hat{g}(t):=
\begin{cases}
&g(x,t)=\chi_{\Omega}(x)f(t)+(1-\chi_{\Omega}(x)) \hat{f}(t)\; \ \ \text{if}\; t> 0,\\
&0\; \ \ \hspace{6.1cm} \text{if}\; t\leq 0.
\end{cases}
\end{align*}

To study problem \eqref{pr3}, we introduce the Euler-Lagrange functional
$\mathcal {I}_{\varepsilon}: \mathcal {W}_{\varepsilon}\to \mathbb R$  by
$$ \mathcal {I}_{\varepsilon}(u)=\sum_{\wp\in \{p,q\}}\frac{1}{\wp}\|u\|_{W_{Z,\varepsilon}^{s,\wp}(\mathbb R^N)}^{\wp}-\frac{1}{2}\int\limits_{\mathbb R^N}\int\limits_{\mathbb R^N}\frac{G(u(y))G(u(x))}{|x-y|^\mu}dy dx.$$
 By conditions $(\mathcal {Z}_1), (\mathcal {Z}_2)$ and $(f_1),$
 the functional $\mathcal {I}_{\varepsilon}$ is well defined on $\mathcal {W}_{\varepsilon}$
 and of class $C^{2}(\mathcal {W}_{\varepsilon})$. Moreover, the critical points
 of $\mathcal {I}_{\varepsilon}$ are exactly the (weak) solutions of problem \eqref{pr3}.
 Associated to
the energy functional $\mathcal {I}_{\varepsilon},$ we denote the Nehari
manifold $\mathscr{N}_{\varepsilon}$  by
$$ \mathscr{N}_{\varepsilon}=\{v\in \mathcal {W}_{\varepsilon}\setminus \{0\}\,:
\,\langle \mathcal {I}'_{\varepsilon}(v),v\rangle=0\},$$
where
\begin{align*}
    \langle\mathcal {I}'_{\varepsilon}(v),\varphi\rangle&=\sum_{\wp\in \{p,q\}}\int\limits_{\mathbb R^{2N}}\frac{|v(x)-v(y)|^{\wp-2}(v(x)-v(y))(\varphi(x)-\varphi(y))}{|x-y|^{N+\wp s}}dxdy\\
    &\quad +\int\limits_{\mathbb R^N}Z(\varepsilon x)(|v|^{p-2}v+|v|^{q-2}v)\varphi dx
    -\int\limits_{\mathbb R^N}\int\limits_{\mathbb R^N}\frac{G(v(y))g(v(x))\varphi(x)}{|x-y|^\mu}dy dx
\end{align*}
for every $v,\varphi\in \mathcal {W}_{\varepsilon}.$
\begin{lemma}\label{lem6a}
Suppose that conditions $(f_1)-(f_5)$ hold. The the following statements hold:
\begin{itemize}
\item[$(g_1)$] $\mathop {\lim }\limits_{t \to 0^+}\frac{g(x, t)}{t^{q-1}} = 0$
uniformly with respect to $x\in\mathbb{R}^N$;
\item[$(g_2)$] \begin{itemize}\item[$(i)$] $0 < g(x,t)\le f(t)$ for every $t>0,$ for every $x\in\mathbb{R}^N;$
\item[$(ii)$] $g(x,t)=0$ for every $t\leq0,$ for every $x\in\mathbb{R}^N$;
\end{itemize}
\item[$(g_3)$]
\begin{itemize}
\item[$(i)$] $0 < \theta G(x,t) \leq g(x,t)t,$ for every $x\in {\Omega}$ and $t>
0$;
\item[$(ii)$] $0\leq q G(x,t)\leq g(x,t)t\leq \frac{Z_0}{\hbar_0}(t^{p}+t^{q})$ for every $x\in \mathbb{R}^N\setminus{\Omega}$ and $t>
0$;
\end{itemize}
\item[$(g_4)$]  for every $x\in \Omega $, the function $t \mapsto \frac{g(x,t)}{t^{\frac{q}{2}- 1}}$
is strictly increasing on $(0, +\infty)$;
\item[$(g_5)$] for every $x\in\mathbb{R}^N\setminus{\Omega}$, function $t \mapsto \frac{g(x,t)}{t^{\frac{q}{2}- 1}}$ is strictly increasing on $(0,a)$.
\end{itemize}
\end{lemma}
\begin{proof}
We shall only give the proof of $(g_3)-(ii).$ The rest of the properties can be verified by the definition of
$g.$ Using $(f_5),$ we obtain that
$$\frac{f(t)}{t^q}\leq \frac{f(a)}{a^q}= \frac{Z_0}{\hbar_0} \ \ \mbox{for every} \ \ t\in[0,a].$$
Consequently, $$g(x,t)t=f(t)t\leq \frac{Z_0}{\hbar_0}t^q\leq\frac{Z_0}{\hbar_0}(t^p+t^q)\
\hbox{ for every } t\in[0,a].$$
 If $t\in(a,+\infty),$
we have $$g(x,t)=\frac{Z_0}{\hbar_0}t^q<\frac{Z_0}{\hbar_0}(t^p+t^q).$$
From  $(f_2)$ and $(g_2),$  we obtain that
$$g(x,t)t=f(t)t\geq \theta F(t)> qF(t)\geq qG(x,t)>0\
\hbox{ for every } t\in[0,a].$$
In addition, if $t\in (a,+\infty),$ we have
$$g(x,t)=\hat{f}(t)=\frac{Z_0}{\hbar_0}t^{q-1}.$$
Hence
$$G(x,t)=\frac{Z_0}{q\hbar_0}t^q$$
and
$$g(x,t)t=\frac{Z_0}{\hbar_0}t^q=qG(x,t).$$
This completes the proof of Lemma \ref{lem6a}.
\end{proof}

\begin{lemma}\label{lmn1}(see \cite[Proposition 3.1]{Liang1})
Suppose that conditions $(\mathcal {Z}_1)-\mathcal {Z}_2$ and
$(f_1)-(f_3)$ hold. Then there is a real number $\mathfrak r_{*}>0$ such that
    $$\|u\|_{\mathcal {W}_{\varepsilon}}\ge \mathfrak r_{*}>0\
    \mbox{for every } u\in \mathscr{N}_{\varepsilon}.$$
\end{lemma}
\begin{lemma}\label{lmm1}
Suppose that conditions $(\mathcal {Z}_1)-(\mathcal {Z}_2)$ and $(f_1)-(f_5)$ hold. Then
$\mathcal {I}_{\varepsilon}$ satisfies the following the geometric conditions:
\begin{itemize}
\item[$(i)$] There exist  real numbers $\alpha_*>0,\rho_*>0$ such that for every $u\in  \mathcal {W}_{\varepsilon}:$ $\|u\|_{\mathcal {W}_{\varepsilon}}=\rho_*,$ we have $\mathcal {I}_{\varepsilon}(u)\ge \alpha_*>0;$
\item[$(ii)$]  There exists $u\in \mathcal {W}_{\varepsilon}$ such that $\|u\|_{\mathcal {W}_{\varepsilon}}>\rho_*$ and $\mathcal {I}_{\varepsilon}(u)<0.$
\end{itemize}
\end{lemma}
\begin{proof}
One can apply a similar discussion as in Lemmas~\ref{lm2}
and~\ref{lm3}, combined with the Trudinger-Moser inequality, so we shall omit the details here.
\end{proof}
By virtue of Lemma \ref{lmn1} and the Mountain Pass Theorem, there
exists a $(PS)_{c_\varepsilon}$ sequence $\{\mathfrak
u_{n}\}_n\subset \mathcal {W}_{\varepsilon},$ that is,
$$\mathcal {I}_\varepsilon (\mathfrak u_n)\rightarrow c_\varepsilon
\ \ \mbox{and} \ \ \mathcal {I}'_\varepsilon (\mathfrak u_n)\rightarrow 0,$$ where
$$c_\varepsilon : = \inf_{\gamma  \in \Gamma }\max_{t \in [0,1]} {\mathcal {I}_\varepsilon }(\gamma (t)) $$
and $\Gamma = \{\gamma  \in (C^0[0,1],\mathcal {W}_\varepsilon(\mathbb R^N))
:\gamma(0) = 0, \mathcal {I}_\varepsilon (\gamma (1))\} < 0$.

We give the definition of  $\mathcal {T}_\varepsilon^+$ as follows:
$$\mathcal {T}_\varepsilon^+:= \{u\in \mathcal {W}_\varepsilon(\mathbb{R}^N):|\text{supp}(u^+)\cap \Omega_\varepsilon|>0\}
\subset \mathcal {W}_\varepsilon(\mathbb{R}^N),$$
$\Omega_\varepsilon:= \{x\in\mathbb{R}^N:\varepsilon x\in\Omega\}.$
Let $\mathbb{S}_\varepsilon$ be a the unit sphere in $\mathcal {W}_\varepsilon(\mathbb{R}^N)$
and  denote by $\mathbb{S}_\varepsilon^+=\mathbb{S}_\varepsilon\cap \mathcal {T}_\varepsilon^+.$
Note that
$\mathbb{S}_\varepsilon^+$ is an incomplete $C^{1,1}$-manifold of
codimension 1, modelled on $\mathcal {W}_\varepsilon(\mathbb{R}^N)$ and contained
in the open $\mathcal {T}_\varepsilon^+$. Thus, $\mathcal {W}_\varepsilon=T_u\mathbb{S}_\varepsilon^+\oplus \mathbb{R}u,$
for every $u\in\mathbb{S}_\varepsilon^+,$
where
\begin{align*}
T_u\mathbb{S}_{\varepsilon}^+ & =\{v\in \mathcal {W}:
\sum_{\wp\in \{p,q\}}\iint\limits_{\mathbb R^{2N}}\frac{|u(x)-u(y)|^{\wp-2}(u(x)-u(y))(\varphi(x)-\varphi(y))}{|x-y|^{N+\wp s}}dxdy \\
& \quad \quad +\int\limits_{\mathbb R^N}Z(\varepsilon x)(|u|^{p-2}u+|u|^{q-2}u)\varphi dx=0\}.
\end{align*}

As in Lemmas \ref{lm7} and \ref{lm7.1}, the following results can be obtained. 
\begin{lemma}\label{lm8.1}
Suppose  that conditions $(\mathcal {Z}_1)-(\mathcal {Z}_2)$ and $(f_1)-(f_5)$ hold. Then the following statements are true:
\begin{itemize}
\item[$(i)$] There exists a unique
$t_uu\in \mathscr{N}_{\varepsilon}$ and $\mathcal {I}_{\varepsilon}(t_uu)=\max_{t>0}\mathcal {I}_{\varepsilon}(tu)$ for every $u\in \mathcal {T}_\varepsilon^+.$  Moreover,
we have $c_{\varepsilon}\geq \zeta>0$ and
$$c_\varepsilon=\inf_{u\in\mathscr{N}_\varepsilon}\mathcal {I}_{\varepsilon}(u)=\inf_{u\in\mathcal {T}_\varepsilon^+}\max_{t>0}\mathcal {I}_\varepsilon(tu)
=\inf_{u\in\mathbb{S}_\varepsilon^+}\max_{t>0}\mathcal {I}_\varepsilon(tu).$$

\item[$(ii)$] $\mathscr{N}_{\varepsilon}$ is bounded away from 0, and there exists $\alpha>0$ such that
$t_u\geq \alpha,$ for every $u\in \mathbb{S}_{\varepsilon}.$
Moreover, for every compact subset $\mathcal {K}\subset\mathbb{S}_{\varepsilon}^+,$ there exists
$C_\mathcal {K}>0$ such that $t_u\leq C_\mathcal {K},$ for every $u\in \mathcal {K}.$

\item[$(iii)$] The continuous map $\hat{m}_{\varepsilon}: \mathcal {W}_{\varepsilon}\rightarrow\mathscr{N}_{\varepsilon}$ is
given by $\hat{m}_{\varepsilon}(u)=t_uu$ and $m_{\varepsilon}:=\hat{m}_{\varepsilon}|_{\mathbb{S}_\varepsilon^+}$ is a homeomorphism between $\mathbb{S}_{\varepsilon}^+$ and $\mathscr{N}_{\varepsilon},$
and $m_{\varepsilon}^{-1}(u)=\frac{u}{\|u\|_{\mathcal {W}_\varepsilon}}.$
\end{itemize}
\end{lemma}
Exploring the functional $\Theta_\varepsilon(u):=\mathcal
{I}_\varepsilon(m_\varepsilon(u))$, together with argument similar  to the proof of   Szulkin and Weth~\cite[Corollary 3.3]{LQ}, we
obtain the following results.
\begin{lemma}\label{lm8.2}
Suppose that conditions $(\mathcal {Z}_1)-(\mathcal {Z}_2)$ and $(f_1)-(f_5)$ hold. Then the following statements are true:
\begin{itemize}
\item[$(i)$] If $\{\mathfrak u_n\}_n$ is a $(PS)_{c_\varepsilon}$ sequence for $\Theta_{\varepsilon},$ then $\{m_{\varepsilon}(\mathfrak u_n)\}_n$ is
a $(PS)_{c_\varepsilon}$ sequence for $\mathcal {I}_{\varepsilon}.$ Moreover, if $\{\mathfrak u_n\}_n\subset\mathscr{N}_{\varepsilon}$ is a bounded $(PS)_{c_\varepsilon}$ sequence for $\mathcal {I}_{\varepsilon},$
then $\{m_{\varepsilon}^{-1}(\mathfrak u_n)\}_n$ is a $(PS)_{c_\varepsilon}$ sequence for $\Theta_{\varepsilon}.$

\item[$(ii)$]  $u$ is a critical point of $\Theta_{\varepsilon}$ if and only if $m_{\varepsilon}(u)$ is a nontrivial critical point of $\mathcal {I}_{\varepsilon}.$ Furthermore, $\inf_{\mathscr{N}_{\varepsilon}}\mathcal {I}_{\varepsilon}=\inf_{\mathbb{S}_{\varepsilon}^+}\Theta_{\varepsilon}.$
\end{itemize}
\end{lemma}

\begin{lemma}\label{lm10b}
 $c_\varepsilon$ and $c_{\vartheta}$ satisfy the following inequalities
\begin{equation}\label{e3.6a}
\limsup_{\varepsilon\rightarrow0^+}c_\varepsilon\leq c_{\vartheta}\leq C_{\gamma_1}.
\end{equation}
\end{lemma}
\begin{proof}
Let $\varphi\in C_0^{\infty}(\mathbb R^N)$ such that $\varphi\equiv 1$ on
$B_{\delta/2}(0)$, supp$(\varphi)\subset B_\delta(0)\subset \Omega$
for some $\delta>0$ and $\varphi\equiv0$ on $B_\delta(0)^c.$ We define
$$u_\varepsilon(x):= \varphi(\varepsilon x)u(x) \ \ \text{for every} \ \  \varepsilon>0,$$
where $u$ is the ground state solution of problem \eqref{ct1a}
obtained in Proposition \ref{pro1}.
We know that
supp$(u_\varepsilon)\subset \Omega_\varepsilon$ and
$u_\varepsilon\rightarrow u$ in $\mathcal {W}$ (see
Ambrosio and Isernia~\cite[Lemma 2.4]{Amb}). If we  assume that  $t_\varepsilon>0$
with $t_\varepsilon u_\varepsilon\in \mathscr{N}_\varepsilon,$ then
\begin{align*}
c_\varepsilon  \leq \mathcal {I}_\varepsilon(t_\varepsilon v_\varepsilon)
&=\frac{t_\varepsilon^p}{p}\|\mathfrak u_n\|_{W_{Z,\varepsilon}^{s,p}(\mathbb R^N)}^{p}
+\frac{t_\varepsilon^q}{q}\|\mathfrak u_n\|_{W_{Z,\varepsilon}^{s,q}(\mathbb R^N)}^{q}
-\frac{1}{2}\int\limits_{\mathbb R^N}\int\limits_{\mathbb R^N}\frac{G(t_\varepsilon u_\varepsilon(y))G(t_\varepsilon u_\varepsilon(x))}{|x-y|^\mu}dy dx \\
& =\frac{t_\varepsilon^p}{p}\|\mathfrak u_n\|_{W_{Z,\varepsilon}^{s,p}(\mathbb R^N)}^{p}
+\frac{t_\varepsilon^q}{q}\|\mathfrak u_n\|_{W_{Z,\varepsilon}^{s,q}(\mathbb R^N)}^{q}
-\frac{1}{2}\int\limits_{\mathbb R^N}\int\limits_{\mathbb R^N}\frac{F(t_\varepsilon u_\varepsilon(y))F(t_\varepsilon u_\varepsilon(x))}{|x-y|^\mu}dy dx.
\end{align*}
Since $t_\varepsilon u_\varepsilon\subset\mathscr{N}_\varepsilon,$ we have
\begin{equation}\label{e3.7}
t_\varepsilon^p\|\mathfrak u_n\|_{W_{Z,\varepsilon}^{s,p}(\mathbb R^N)}^{p}
+t_\varepsilon^q\|\mathfrak u_n\|_{W_{Z,\varepsilon}^{s,q}(\mathbb R^N)}^{q}
=\int\limits_{\mathbb R^N}\int\limits_{\mathbb R^N}\frac{F(t_\varepsilon u_\varepsilon(y))f(t_\varepsilon u_\varepsilon(x))t_\varepsilon u_\varepsilon(x)}{|x-y|^\mu}dy dx.
\end{equation}
Moreover, it follows from $(f_3)$ that
\begin{align}\label{e3.8}
 \mathcal {I}_{\varepsilon}(t_\varepsilon v_\varepsilon)
& \nonumber=\frac{t_\varepsilon^p}{p}\|\mathfrak u_n\|_{W_{Z,\varepsilon}^{s,p}(\mathbb R^N)}^{p}
+\frac{t_\varepsilon^q}{q}\|\mathfrak u_n\|_{W_{Z,\varepsilon}^{s,q}(\mathbb R^N)}^{q}
-\frac{1}{2}\int\limits_{\mathbb R^N}\int\limits_{\mathbb R^N}\frac{F(t_\varepsilon u_\varepsilon(y))F(t_\varepsilon u_\varepsilon(x))}{|x-y|^\mu}dy dx \\
& \geq\int_{\mathbb{R}^N}\int_{\mathbb{R}^N}\frac{F(t_\varepsilon u_\varepsilon(y)}{|x-y|^\mu}
\Big[\frac{1}{q}f(t_\varepsilon u_\varepsilon)t_\varepsilon
u_\varepsilon dx-\frac{1}{2}F(t_\varepsilon u_\varepsilon)\Big]dx \geq0.
\end{align}
This fact implies that $\{t_\varepsilon\}$ is bounded, as $\varepsilon\rightarrow0^+.$
Indeed, suppose to the contrary, that $\{t_\varepsilon\}$ is unbounded when $\varepsilon\rightarrow0^+.$ Consequently,  together with $(f_4),$ we would have
\begin{equation*}
\mathcal {I}_\varepsilon(t_\varepsilon u_\varepsilon)
\geq  \frac{t_\varepsilon^p}{p}\|u_\varepsilon\|_{\vartheta,W^{s,p}_{Z,\varepsilon}(\mathbb{R}^N)}^p
+\frac{t_\varepsilon^q}{q}\|u_\varepsilon\|_{W_{Z,\varepsilon}^{s,q}(\mathbb R^N)}^{q}
-\gamma_1^2t_\varepsilon^{2\theta}\||u_\varepsilon|^\theta\|_{L^{\frac{2N}{2N-\mu}}}^2\rightarrow-\infty
\end{equation*}
which is a contradiction with \eqref{e3.8}.
Thus, we may suppose that $t_\varepsilon\rightarrow t_0,$ as
$\varepsilon\rightarrow0^+.$  Using the Vitali's Theorem, we obtain
\begin{align*}
\limsup_{\varepsilon\rightarrow0^+} c_\varepsilon
& \leq \frac{t_0^p}{p}\|u\|_{W_{Z,\varepsilon}^{s,p}(\mathbb R^N)}^{p}
+\frac{t_0^q}{q}\|u\|_{W_{Z,\varepsilon}^{s,q}(\mathbb R^N)}^{q}
-\frac{1}{2}\int\limits_{\mathbb R^N}\int\limits_{\mathbb R^N}\frac{F(t_0 u(y))
F(t_0 u(x))}{|x-y|^\mu}dy dx
 =\mathcal {E}_{\vartheta}(t_0u).
\end{align*}
We shall now verify that $t_0=1.$
For every $\varepsilon>0,$ there exists $t_\varepsilon>0$ such that
$$\mathcal {I}_\varepsilon(t_\varepsilon u_\varepsilon)=\max_{t\geq0}\mathcal {I}_\varepsilon(tu_\varepsilon).$$
Therefore, $\langle\mathcal {I}'_\varepsilon(t_\varepsilon u_\varepsilon),u_\varepsilon\rangle=0$
and we have
$$t_\varepsilon^p\|u_\varepsilon\|_{W_{Z,\varepsilon}^{s,p}(\mathbb R^N)}^{p}
+t_\varepsilon^q\|u_\varepsilon\|_{W_{Z,\varepsilon}^{s,q}(\mathbb R^N)}^{q}
=\int\limits_{\mathbb R^N}\int\limits_{\mathbb R^N}\frac{F(t_\varepsilon u_\varepsilon(y))f(t_\varepsilon u_\varepsilon(x))t_\varepsilon u_\varepsilon(x)}{|x-y|^\mu}dy dx$$ which means that
\begin{equation}\label{e3.9}
\|u_\varepsilon\|_{W_{Z,\varepsilon}^{s,p}(\mathbb R^N)}^{p}
+t_\varepsilon^{q-p}\|u_\varepsilon\|_{W_{Z,\varepsilon}^{s,q}(\mathbb R^N)}^{q}
=\int\limits_{\mathbb R^N}\int\limits_{\mathbb R^N}\frac{F(t_\varepsilon u_\varepsilon(y))f(t_\varepsilon u_\varepsilon(x))u_\varepsilon(x)}{t_\varepsilon^{p-1}|x-y|^\mu}dy dx.
\end{equation}
Passing to the limit as $\varepsilon\rightarrow0$ in \eqref{e3.9} and
using the fact that $u_\varepsilon\rightarrow u$ in
$W^{s,\frac{N}{s}}(\mathbb{R}^N)$
(see  Ambrosio and Isernia~\cite[Lemma
2.4]{Amb}), we obtain
\begin{equation*}
\|u\|_{W_{Z,\varepsilon}^{s,p}(\mathbb R^N)}^{p}
+t_0^{q-p}\|u\|_{W_{Z,\varepsilon}^{s,q}(\mathbb R^N)}^{q}
=\int\limits_{\mathbb R^N}\int\limits_{\mathbb R^N}\frac{F(t_0 u(y))f(t_0 u(x))u(x)}{t_0^{p-1}|x-y|^\mu}dy dx.
\end{equation*}
Together with this fact, $u\in\mathscr{M}_{\vartheta}$ and $(f_5),$ we deduce that $t_0 = 1.$
Therefore,
$$\limsup_{\varepsilon\rightarrow 0^+}c_\varepsilon\leq
\mathcal {E}_{\vartheta}(v)=c_{\vartheta}.$$ Combining with Lemma \ref{lm3a}, we can obtain that the last
inequality hold in \eqref{e3.6a}. This completes the proof of Lemma
\ref{lm10b}.
\end{proof}

\begin{lemma}\label{lem7}
The functional  $\mathcal {I}_\varepsilon$ satisfies the Palais-Smale
condition at level $c_\varepsilon,$ for every $\varepsilon\in(0,\varepsilon_0).$
\end{lemma}

\begin{proof}
Let $\{\mathfrak u_n\}_n\subset \mathcal {W}_\varepsilon(\mathbb{R}^N)$ be a $(PS)_{c_\varepsilon}$ sequence for $\mathcal {I}_\varepsilon,$
i.e.,
\begin{equation*}
\mathcal {I}_\varepsilon(\mathfrak u_n)\rightarrow c_\varepsilon
\ \ \mbox{and} \ \ \mathcal {I}'_\varepsilon(\mathfrak u_n)\rightarrow0 
\ \ \mbox{as} \ \ n\rightarrow\infty.
\end{equation*}
We shall complete the proof of this lemma by the following two claims.

\medskip

$\textbf{Claim 1.}$
$\{\mathfrak u_n\}_n$ is bounded in $\mathcal {W}_\varepsilon(\mathbb R^N).$

Indeed, by $(f_3),$  we have
\begin{align*}
c_\varepsilon+o_n(1)& = \mathcal {I}_\varepsilon(\mathfrak u_n)-\frac{1}{2q}\langle\mathcal {I}'_\varepsilon(\mathfrak u_n),\mathfrak u_n\rangle
 =(\frac{1}{p}-\frac{1}{2q})\|\mathfrak u_n\|_{W_{Z,\varepsilon}^{s,p}(\mathbb R^N)}^{p}
+\frac{1}{2q}\|\mathfrak u_n\|_{W_{Z,\varepsilon}^{s,q}(\mathbb R^N)}^{q} \\
& \quad + \frac{1}{2}\int_{\mathbb{R}^N}\int_{\mathbb{R}^N}
\frac{F(\mathfrak u_n(y))}{|x-y|^{\mu}}[\frac{1}{q}f(\mathfrak u_n(x))\mathfrak u_n(x)-F(\mathfrak u_n(x))]dydx \\
& \geq\frac{1}{2q}(\|\mathfrak u_n\|_{W_{Z,\varepsilon}^{s,p}(\mathbb R^N)}^{p}
+\|\mathfrak u_n\|_{W_{Z,\varepsilon}^{s,q}(\mathbb R^N)}^{q}) =\frac{1}{2q}\|\mathfrak u_n\|_{\mathcal {W}_\varepsilon}
\end{align*}
which implies that $\{\mathfrak u_n\}_n$ is bounded in $\mathcal {W}_\varepsilon(\mathbb R^N).$ Moreover,
we obtain
\begin{equation}\label{e3.10}
\limsup_{n\rightarrow\infty}\|\mathfrak u_n\|_{\mathcal {W}_\varepsilon}^p
\leq 2qc_\varepsilon.
\end{equation}
Together with this fact and Lemma \ref{lm3a}, we obtain that
\begin{align}\label{e3.10a}
 \limsup_{n\rightarrow\infty}\|\mathfrak u_n\|_{\mathcal {W}_\varepsilon}
\leq 2qc_{\vartheta}
\leq 2qC_{\gamma_1}
&= 2q a\left(1-\frac{N}{2\theta s}\right)
    \left(\frac{a N}{2\theta sb}\right)^{N/(2\theta s-N)}:=\mathcal {G}
\end{align}
    for $\gamma_1\geq \max\{\gamma^\ast,a\},$ where $\gamma^\ast, a$ are given in \eqref{1a} and Lemma
\ref{lm3a}, respectively.  Therefore,  going to a subsequence if necessary, we may assume that
$\mathfrak u_n\rightharpoonup u$ in $\mathcal {W}_\varepsilon,$
$\mathfrak u_n\rightarrow u$ in $L_{loc}^q(\mathbb{R}^N)$ for every
$q\in[\frac{N}{s},+\infty)$ and $\mathfrak u_n \rightarrow u$ a.e. in
$\mathbb{R}^N.$

\medskip

\textbf{Claim 2.} $(PS)_{c_\varepsilon}$  condition  holds in $\mathcal {W}_\varepsilon.$

We  shall divide the proof into three steps.
We first verify that $u$ is a critical point of $\mathcal {I}_\varepsilon.$
For every $\varphi  \in C_c^\infty \left( {{\mathbb{R}^N}} \right),$ we have
\begin{align*}
& \sum_{\wp\in \{p,q\}}\iint\limits_{\mathbb R^{2N}}\frac{|\mathfrak u_n(x)-\mathfrak u_n(y)|^{\wp-2}(\mathfrak u_n(x)-\mathfrak u_n(y))(\varphi(x)-\varphi(y))}{|x-y|^{N+\wp s}}dxdy \\
& \rightarrow
\sum_{\wp\in \{p,q\}}\iint\limits_{\mathbb R^{2N}}\frac{|u(x)-u(y)|^{\wp-2}(u(x)-u(y))(\varphi(x)-\varphi(y))}{|x-y|^{N+\wp s}}dxdy
\end{align*}
and $$\int_{\mathbb{R}^N}Z(\varepsilon x)|\mathfrak u_n|^{\wp-2}\mathfrak u_n\varphi dx\rightarrow
\int_{\mathbb{R}^N}Z(\varepsilon x)|u|^{\wp-2}u\varphi dx 
\ \ \mbox{for every} \ \ \wp\in \{p,q\}.$$
\medskip

$\textbf{Step 1.}$ We prove that
\begin{equation}\label{4c}
\int_{\mathbb{R}^N}\int_{\mathbb{R}^N}\frac{G(\mathfrak u_n(y))g(\mathfrak u_n(x))\varphi(x)}{|x-y|^{\mu}}dydx
\rightarrow \int_{\mathbb{R}^N}\int_{\mathbb{R}^N}\frac{G(u(y))g(u(x))\varphi(x)}{|x-y|^{\mu}}dydx.
\end{equation}
Since the boundedness of $\{G(\varepsilon x,\mathfrak u_n)\}_n$ in ${L^{\frac{{2N}}{{2N - \mu }}}}\left( {{\mathbb{R}^N}} \right),$ ${\mathfrak u_n} \to u$ a.e. in ${{\mathbb{R}^N},}$ and $t\mapsto G\left( { \cdot ,t} \right)$ is continuous, hence
\begin{align*}
G\left( {\varepsilon x,{\mathfrak u_n}} \right) \rightharpoonup G\left( {\varepsilon x,u} \right) \;\text{in}\; {L^{\frac{{2N}}{{2N - \mu }}}}\left( {{\mathbb{R}^N}} \right). \ \
\end{align*}
From Lemma \ref{lm2n}, we obtain the linear bounded operator
\begin{align*}
\frac{1}{{{{\left| x \right|}^\mu }}} * F \in {L^{\frac{{2N}}{\mu }}}\left( {{\mathbb{R}^N}} \right)  \ \ \text{for every}\;F \in {L^{\frac{{2N}}{{2N - \mu }}}}\left( {{\mathbb{R}^N}} \right),
\end{align*}
from ${L^{\frac{{2N}}{{2N - \mu }}}}\left( {{\mathbb{R}^N}} \right)$ to ${L^{\frac{{2N}}{\mu }}}\left( {{\mathbb{R}^N}} \right).$ Therefore,
$\frac{1}{{{{\left| x \right|}^\mu }}}\ast G\left( {\varepsilon y,{\mathfrak u_n}} \right) \rightharpoonup \frac{1}{{{{\left| x \right|}^\mu }}}\ast G\left( {\varepsilon y,u} \right)$ in ${L^{\frac{{2N}}{\mu }}}\left( {{\mathbb{R}^N}} \right).$
Since $g\left( {{\mathfrak u_n}} \right) \to g\left( u \right)$ in $L_{loc}^\nu\left( {{\mathbb{R}^N}} \right),$ for every $\nu \in [p, +\infty),$ we deduce that \eqref{4c} is true.
Consequently, in view of $\langle{{\mathcal {I}'_\varepsilon}(\mathfrak u_{n}),\phi}\rangle  = {o_n}(1),$
for every  $\phi\in C_c^\infty ({{\mathbb{R}^N}}),$ we obtain $\langle{\mathcal {I}'_\varepsilon }(u),\phi\rangle = 0,$ for every $\phi\in C_c^\infty (\mathbb{R}^N).$ Since  $C_c^\infty (\mathbb{R}^N)$
is dense in ${\mathcal {W}_\varepsilon },$ we obtain that $u$ is a critical point of $\mathcal {I}_\varepsilon.$

\medskip

$\textbf{Step 2.}$ We shall prove that for every $\xi,$ there exists $R = R(\xi)>0$ such that
\begin{equation}\label{e3.11}
\limsup_{n\rightarrow\infty}\int_{B_R^c}
\Big(\sum_{\wp\in \{p,q\}}\int\limits_{\mathbb R^{N}}\frac{|u(x)-u(y)|^{\wp}}{|x-y|^{N+\wp s}}dy
+Z(\varepsilon x)|\mathfrak u_n|^{p}+Z(\varepsilon x)|\mathfrak u_n|^{q}\Big)dx<\xi.
\end{equation}
For every
$R > 0,$ let ${\phi _R} \in C_c^\infty \left(
{{\mathbb{R}^N}} \right)$ such that $0 \leq {\phi _R} \leq
1,{\phi _R} = 0$ in $B_R(0),{\phi _R} = 1$ in $B_{2R}^c( 0),$ and
$\left| {\nabla {\phi _R}} \right| \leq \frac{C}{R}$ for some
constant $C>0$ independent of $R.$ Since $\{\phi_R\mathfrak u_n\}_n$ is bounded
in ${\mathcal {W}_\varepsilon}(\mathbb R^N),$ it follows that $\langle \mathcal
{I}'_\varepsilon (\mathfrak u_n), \phi_R\mathfrak u_n\rangle = {o_n}(1),$ i.e.,
\begin{align}\label{p10}
&\nonumber\int_{\mathbb{R}^N}\int_{\mathbb{R}^N}\frac{|\mathfrak u_n(x)-\mathfrak u_n(y)|^p}{|x-y|^{N+sp}}\phi_R(x)dy dx
+\int_{\mathbb{R}^N}\int_{\mathbb{R}^N}\frac{|\mathfrak u_n(x)-\mathfrak u_n(y)|^q}{|x-y|^{N+sq}}\phi_R(x)dy dx \\
&\nonumber \quad + \int_{\mathbb{R}^N}Z(\varepsilon x)|\mathfrak u_n(x)|^pdx
+\int_{\mathbb{R}^N}Z(\varepsilon x)|\mathfrak u_n(x)|^qdx\\
& \nonumber= o_n(1)+\int_{\mathbb{R}^N}\int_{\mathbb{R}^N}\frac{G(\mathfrak u_n(y))g(\mathfrak u_n(x))\phi_R(x)\mathfrak u_n(x)}{|x-y|^{\mu}}dydx \\
&\nonumber\quad -
\iint_{\mathbb{R}^{2N}}\frac{|\mathfrak u_n(x)-\mathfrak u_n(y)|^{p-2}(\mathfrak u_n(x)-\mathfrak u_n(y))(\phi_R(x)-\phi_R(y))}{|x-y|^{N+sp}} \mathfrak u_n(y)dxdy \\
&\quad -
\iint_{\mathbb{R}^{2N}}\frac{|\mathfrak u_n(x)-\mathfrak u_n(y)|^{q-2}(\mathfrak u_n(x)-\mathfrak u_n(y))(\phi_R(x)-\phi_R(y))}{|x-y|^{N+sq}}\mathfrak u_n(y)dxdy.
\end{align}
By Lemma \ref{lma}, then there exists $\hbar_0>0$ such that
\begin{equation}\label{4b}
\frac{\sup_{\mathfrak u_n\in W^{s,p}(\mathbb{R}^N)}\Big|{|x|^{-\mu}}*G(\mathfrak u_n)\Big|}{\hbar_0}<\frac{1}{2}.
\end{equation}
Let $R>0$ be such that $\Omega_\varepsilon\subset B_R.$
By the definition of $\phi_R, (g_3)-(ii)$ and \eqref{4b}, we get
\begin{align}\label{p12}
& \nonumber\int_{\mathbb{R}^N}\int_{\mathbb{R}^N}\frac{|\mathfrak u_n(x)-\mathfrak u_n(y)|^p}{|x-y|^{N+sp}}\phi_R(x)dy dx
+\int_{\mathbb{R}^N}\int_{\mathbb{R}^N}\frac{|\mathfrak u_n(x)-\mathfrak u_n(y)|^q}{|x-y|^{N+sq}}\phi_R(x)dy dx \\
&\nonumber \quad + \frac{1}{2}\int_{\mathbb{R}^N}Z(\varepsilon x)(|\mathfrak u_n(x)|^p+|\mathfrak u_n(x)|^q)dx  \\
&\nonumber\leq{o_n}-
\iint_{\mathbb{R}^{2N}}\frac{|\mathfrak u_n(x)-\mathfrak u_n(y)|^{p-2}(\mathfrak u_n(x)-\mathfrak u_n(y))(\phi_R(x)-\phi_R(y))}{|x-y|^{N+sp}} \mathfrak u_n(y)dxdy \\
&\quad -
\iint_{\mathbb{R}^{2N}}\frac{|\mathfrak u_n(x)-\mathfrak u_n(y)|^{q-2}(\mathfrak u_n(x)-\mathfrak u_n(y))(\phi_R(x)-\phi_R(y))}{|x-y|^{N+sq}}\mathfrak u_n(y)dxdy.
\end{align}
For $\wp\in \{p,q\},$ by virtue of  the H\"{o}lder
inequality and the boundedness of $\{\mathfrak u_n\}_n$ in
$\mathcal {W}_\varepsilon,$ we have
\begin{align}\label{p12a}
& \nonumber\iint_{\mathbb{R}^{2N}}\frac{|\mathfrak u_n(x)-\mathfrak u_n(y)|^{\wp-2}(\mathfrak u_n(x)-\mathfrak u_n(y))
(\phi_R(x)-\phi_R(y))}{|x-y|^{N+\wp s}}\mathfrak u_n(y)dxdy\\
&\leq C\Big(\iint_{\mathbb{R}^{2N}}\frac{|\phi_R(x)-\psi_R(y)|^{\wp}}{|x-y|^{N+\wp s}}
|\mathfrak u_n(y)|^{\wp}dy\Big)^{\frac{1}{\wp}}.
\end{align}
Next, by the definition of $\phi_R,$ polar
coordinates and the boundedness of $\{\mathfrak u_n\}_n$ in
$\mathcal {W}_\varepsilon,$ we obtain
\begin{align}\label{p13}
&\nonumber\iint_{{\mathbb{R}^{2N}}} {\frac{{{{\left| {{\phi _R}\left( x \right) - {\phi _R}\left( y \right)} \right|}^{\wp}}}}{{{{\left| {x - y} \right|}^{N+\wp s}}}}{{\left| {{\mathfrak u_n}\left( x \right)} \right|}^{\wp}}dxdy}  \\
&\nonumber= \int_{{\mathbb{R}^N}} {\int_{\left| {y - x} \right| > R} {\frac{{{{\left| {{\phi_R}\left( x \right) - {\phi_R}\left( y \right)} \right|}^{\wp}}}}{{{{\left| {x - y} \right|}^{N+\wp s}}}}{{\left| {{\mathfrak u_n}
\left( x \right)} \right|}^{\wp}}dxdy  } }
\nonumber\quad +\int_{{\mathbb{R}^N}} {\int_{\left| {y - x} \right| \leq R} {\frac{{{{\left| {{\phi _R}\left( x \right) - {\phi_R}\left( y \right)} \right|}^{\wp}}}}{{{{\left| {x - y} \right|}^{N+\wp s}}}}{{\left| {{\mathfrak u_n}\left( x \right)} \right|}^{\wp}}dxdy} }  \\
&\nonumber\leq C\int_{\mathbb{R}^N} {{{\left| {{\mathfrak u_n}\left( x
\right)} \right|}^{\wp}}\left( {\int_{\left| {y - x} \right| > R}
{\frac{{dy}}{{{{\left| {x - y} \right|}^{N+\wp s}}}}} } \right)dx
+ \frac{C}{{{R^{\wp}}}}} \int_{{\mathbb{R}^N}} {{{\left| {{\mathfrak u_n}\left( x \right)} \right|}^{\wp}}} \left( {\int_{\left| {y - x} \right| \leq R} {\frac{{dy}}{{{{\left| {x - y} \right|}^{N+\wp s - {\wp}}}}}} } \right)dx  \\
&\nonumber\leq C\int_{{\mathbb{R}^N}} {{{\left| {{\mathfrak u_n}\left( x \right)} \right|}^{\wp}}\left( {\int_{\left| z \right| > R} {\frac{{dz}}{{{{\left| z \right|}^{N+\wp s}}}}} } \right)dx + \frac{C}{{{R^{t}}}}} \int_{{\mathbb{R}^N}} {{{\left| {{\mathfrak u_n}\left( x \right)} \right|}^{\wp}}} \left( {\int_{\left| z \right| \leq R} {\frac{{dz}}{{{{\left| z \right|}^{N+\wp s - {\wp}}}}}} } \right)dx  \\
&\nonumber\leq C{\int_{{\mathbb{R}^N}} {\left| {{\mathfrak u_n}\left( x \right)} \right|} ^{\wp}}dx\left( {\int_R^\infty  {\frac{{d\rho }}{{{\rho ^{s\wp+ 1}}}}} } \right) + \frac{C}{{{R^{\wp}}}}\int_{{\mathbb{R}^N}} {{{\left| {{\mathfrak u_n}\left( x \right)} \right|}^{\wp}}} dx\left( {\int_0^R {\frac{{d\rho }}{{{\rho ^{s\wp- {\wp} + 1}}}}} } \right) \\
&\nonumber\leq \frac{C}{{{R^{{s\wp}}}}}{\int_{{\mathbb{R}^N}} {\left| {{\mathfrak u_n}(x)} \right|}^{\wp}}dx + \frac{C}{{{R^{\wp}}}}{R^{ -s\wp + {\wp}}}\int_{{\mathbb{R}^N}} {{{\left| {{\mathfrak u_n}\left( x \right)} \right|}^{\wp}}} dx
  \nonumber\leq \frac{C}{{{R^{s\wp}}}}\int_{{\mathbb{R}^N}} {{{\left| {{\mathfrak u_n}\left( x \right)} \right|}^{\wp}}} dx\\
&\leq \frac{C}{{{R^{s\wp}}}}\rightarrow0 \ \ \mbox{as} \ \ R\rightarrow\infty,
\hbox{ where } \wp\in\{p,q\}.
\end{align}

Gathering \eqref{p12}-\eqref{p13}, we infer that
\eqref{e3.11} is satisfied.

\textbf{Step 3.}
We  verify that ${\mathfrak u_n} \to u$ in $\mathcal {W}_\varepsilon$ as $n \to \infty .$

In view of \eqref{e3.11}, we obtain ${\mathfrak u_n} \to u$ in ${L^\nu}\left( {{\mathbb{R}^N}} \right),$ for every $\nu \in [p, +\infty).$ Indeed, fixed $\xi> 0,$ there exists $R = R(\xi)>0$ such that \eqref{e3.11} holds. Using the compactness embedding $\mathcal {W}_\varepsilon\hookrightarrow\hookrightarrow L^\nu_{loc}(\mathbb{R}^N)$ and  $(\mathcal {Z}_1),$ we see
\begin{align*}
\limsup_{n\rightarrow\infty}\left| {{\mathfrak u_n} - u} \right|_{L^{p}(\mathbb R^N)}^p &\nonumber
= \limsup_{n\rightarrow\infty} \left[ {\left| {{\mathfrak u_n} - u} \right|_{{L^p}\left( {{B_R}\left( 0 \right)} \right)}^p + \left| {{\mathfrak u_n} - u} \right|_{{L^p}\left( {B_R^c\left( 0 \right)} \right)}^p} \right]\\
&\nonumber\leq {2^{p - 1}}\limsup_{n\rightarrow\infty}
\left( {\left| {{\mathfrak u_n}} \right|_{{L^p}\left( {B_R^c\left( 0 \right)} \right)}^p + \left| u \right|_{{L^p}\left( {B_R^c\left( 0 \right)} \right)}^p} \right)\\
&\nonumber\leq\frac{{{2^{p - 1}}}}{{{Z_0}}}\limsup_{n\rightarrow\infty}\int_{B_R^c\left( 0 \right)} {\left( {\int_{{\mathbb{R}^N}} {\frac{{{{\left| {{\mathfrak u_n}\left( x \right) - {\mathfrak u_n}\left( y \right)} \right|}^p}}}{{{{\left| {x - y} \right|}^{2N}}}}dy + Z\left( {\varepsilon x} \right){{\left| {{\mathfrak u_n}} \right|}^p}} } \right)dx}\\
&\nonumber \quad + \frac{{{2^{p - 1}}}}{{{Z_0}}}\int_{B_R^c\left( 0 \right)} {\left( {\int_{{\mathbb{R}^N}} {\frac{{{{\left| {{u}\left( x \right) - {u}\left( y \right)} \right|}^p}}}{{{{\left| {x - y} \right|}^{2N}}}}dy + Z\left( {\varepsilon x} \right){{\left| u \right|}^p}} } \right)dx}
\nonumber< \frac{{{2^p}}}{{{Z_0}}}\xi.
\end{align*}
Due to the arbitrariness of $\xi,$ $\mathfrak u_n\rightarrow u$ in
$L^p(\mathbb{R}^N)$. By interpolation, ${\mathfrak u}_n \to u$ in
${L^\nu}\left( {{\mathbb{R}^N}} \right)$ for every $\nu\in
[p,+\infty),$ as desired. Arguing similarly as in the proof of
\cite[Lemma 22]{Zhang2},
 we can obtain  that
\begin{align}\label{pp19}
\int_{{\mathbb{R}^N}}\int_{{\mathbb{R}^N}}\frac{G(\mathfrak u_n(y))g(\mathfrak u_n(x))}{|x-y|^{\mu}}dydx
\to \int_{{\mathbb{R}^N}}\int_{{\mathbb{R}^N}}\frac{G(u(y))g(u(x))}{|x-y|^{\mu}}dydx.
\end{align}
Therefore, it follows from $\langle\mathcal {I}'_\varepsilon(\mathfrak u_n),\mathfrak u_n\rangle=o_n(1)$
and $\langle\mathcal {I}'_\varepsilon(u),u\rangle=o_n(1)$ that
$$\|\mathfrak u_n\|_{W^{s,p}_{Z,\varepsilon}}^p+\|\mathfrak u_n\|_{W^{s,q}_{Z,\varepsilon}}^q
=
\|u\|_{W^{s,p}_{Z,\varepsilon}}^p+\|u\|_{W^{s,q}_{Z,\varepsilon}}^q+o_n(1).$$
By the Br\'{e}zis-Lieb lemma, we get that
$$\|\mathfrak u_n-u\|_{\mathcal {W}_\varepsilon}^{\wp}=\|\mathfrak u_n\|_{W^{s,\wp}_{Z,\varepsilon}}^{\wp}
-\|u\|_{W^{s,\wp}_{Z,\varepsilon}}^{\wp} +o_n(1) \ \ \mbox{for every} \ \
\wp\in\{p,q\}.$$ Therefore
$$\|\mathfrak u_n-u\|_{W^{s,p}_{Z,\varepsilon}}^p+\|\mathfrak u_n-u\|_{W^{s,q}_{Z,\varepsilon}}^q=o_n(1).$$
Moreover, we obtain $\mathfrak u_n\rightarrow u$ in $\mathcal {W}_\varepsilon.$
This completes the proof of Lemma \ref{lem7}.
\end{proof}

\begin{lemma}\label{lm3.7}(see Ambrosio \cite[Corollary 3.1]{A0})
The functional ${\Phi _\varepsilon}$ satisfies the $(PS)_{c_\varepsilon}$  condition at level $c_\varepsilon,$ for every $\varepsilon\in(0,\varepsilon_0)$ on $\mathbb{S}_{\varepsilon}.$
\end{lemma}

\section{Multiplicity and concentration of positive solutions of problem $\eqref{pr2}$}\label{s5}
\def\theequation{5.\arabic{equation}}
\setcounter{equation}{0}

In this section, we shall prove the main results. To this end, 
we shall give some
notations and useful results which will be used later. Fix $\delta>0$. Let $\mathfrak w$ be
a ground state solution of equation $(\mathcal {Q}_{Z_0})$, so
that $\mathcal {E}_{Z_0}(\mathfrak w)=c_{Z_0}$ and $\mathcal {E}'_{Z_0}(\mathfrak w)=0.$ Let $\eta$ be a
smooth nonincreasing cut-off function in $\mathbb R^+_0$ such that
$\eta(t)=1$ if $0\le t\le {\delta}/{2}$ and $\eta(t)=0$ if $t\ge
\delta.$ For $\varepsilon>0$ and any $y\in \mathscr{M},$ we define
$$\psi_{\varepsilon,y}(x)=\eta(|\varepsilon x-y|)\mathfrak w
\left(\dfrac{\varepsilon x-y}{\varepsilon}\right) ,\quad x\in\mathbb
R^N$$ and $\Phi_{\varepsilon}: \mathscr{M} \to\mathscr{N}_{\varepsilon}$ is given by
$\Phi_{\varepsilon}(y)=t_{\varepsilon}\psi_{\varepsilon,y},$ when
$t_{\varepsilon}>0$ satisfies
$$ \max_{t\ge 0} \mathcal {I}_{\varepsilon}(t\psi_{\varepsilon,y})
= \mathcal {I}_{\varepsilon}(t_{\varepsilon}\psi_{\varepsilon,y}).$$
We obtain that $\Phi_{\varepsilon}(y)$ has compact
 support in $\mathbb R^N$ for every $y\in \mathscr{M}.$

\begin{lemma}\label{lm8}(see Liang et al. \cite[Lemma 5.1]{Liang1})
The function $\Phi_{\varepsilon}$ has the following property
$$ \lim_{\varepsilon\to 0^{+}}\mathcal {I}_{\varepsilon}(\Phi_{\varepsilon}(y))=c_{Z_0},\;\mbox{ uniformly in }\;
y\in \mathscr{M}.$$
\end{lemma}

For any $\delta>0,$ let $\varrho=\varrho(\delta)>0$ be such that $\mathscr{M}_{\delta}\subset B_{\varrho}(0).$ We define the function $\mathscr{X}:\mathbb R^N\to\mathbb R^N$ by
\begin{equation*}
\mathscr{X}(x)= \begin{cases}
x,&\text{if}\; |x|<\varrho,\\
\dfrac{\varrho x}{|x|},&\text{if}\; |x|\ge \varrho.
\end{cases}
\end{equation*}
In what follows, let the barycenter map $\beta_{\varepsilon}: \mathscr{N}_{\varepsilon}\to \mathbb R^N$ be
defined by
$$ \beta_{\varepsilon}(u)=\frac{\int\limits_{\mathbb R^N}\mathscr{X}(\varepsilon x)(|u(x)|^{p}+|u(x)|^{q})dx}{\int\limits_{\mathbb R^N}(|u|^{p}+|u|^{q})dx},
\quad u\in\mathscr{N}_{\varepsilon}.$$

Arguing as in the similar discussion of ~\cite[Lemma~14]{T22},
 we obtain the following result.
\begin{lemma}\label{lm9}
The map $\beta_\varepsilon\circ\Phi_{\varepsilon}$ satisfies the
following limit
\begin{align}\label{nm47}
\lim\limits_{\varepsilon\to 0^{+}}
\beta_{\varepsilon}(\Phi_{\varepsilon}(y))=y,\;\text{uniformly in}\; y\in \mathscr{M}.
\end{align}
\end{lemma}

\begin{lemma}\label{lm10}
 Let $\varepsilon_n\to 0^{+}$ and $\{\mathfrak u_n\}_n\subset \mathscr{N}_{\varepsilon_n}$ satisfy $\mathcal {I}_{\varepsilon_n}(\mathfrak u_n)\to c_{V_0},$
as $n\to\infty$.
 Then there exists a sequence $\{\tilde  y_n\}_n\subset \mathbb R^N$ such that the sequence  $\mathfrak v_n=\mathfrak u_n(\cdot+\tilde  y_n)$ has a subsequence which strongly converges in $\mathcal {W}.$ Furthermore, up to a subsequence,
 $y_n=\varepsilon \tilde y_n\to y\in \mathscr{M}.$
 \end{lemma}
\begin{proof}
Since $\left\langle {{\mathcal {I}_{{\varepsilon _n}}'}(\mathfrak u_n),{\mathfrak u_n}} \right\rangle  = 0$ and ${{\mathcal {I}_{{\varepsilon _n}}}\left( {{\mathfrak u_n}} \right) \to {c_{{Z_0}}}},$ we can see that $\{\mathfrak u_n\}_n$ is a bounded sequence in $\mathcal {W}.$ Indeed, by $(g_3)-(ii)$, we have
\begin{align*}
\mathcal {I}_{\varepsilon_n}(\mathfrak u_n)&=\mathcal {I}_{\varepsilon_n}(\mathfrak u_n)-\frac{1}{2q}\langle \mathcal {I}'_{\varepsilon_n}(\mathfrak u_n),\mathfrak u_n\rangle
 =(\frac{1}{p}-\frac{1}{2q})\|\mathfrak u_n\|_{Z,W^{s,p}_\varepsilon(\mathbb{R}^N)}^p
+\frac{1}{2q}\|\mathfrak u_n\|_{Z,W^{s,q}_\varepsilon(\mathbb{R}^N)}^q  \\
& \quad +\frac{1}{2}\int\limits_{\mathbb R^N}\Big[\frac{1}{|x|^{\mu}}*G(\varepsilon y,\mathfrak u_n(y))\Big]\Big(\frac{1}{q}g(\varepsilon x, u_n(x))\mathfrak u_n(x)-G(\varepsilon x, \mathfrak u_n(x))\Big)dx
\ge \frac{1}{2q}\|\mathfrak u_n\|_{\mathcal {W}_{\varepsilon_n}}.
\end{align*}
Therefore
\begin{align}\label{ctnew3a}
\limsup_{n\to\infty}\|\mathfrak u_n\|_{\mathcal {W}_{\varepsilon}}
\le 2qc_{Z_{0}}.
\end{align}
By conditions $(\mathcal {Z}_1)$ and $(\mathcal {Z}_2),$
we obtain that
\begin{equation}\label{2a}
\|u\|_{W^{s,p}(\mathbb{R}^N)}\leq \min\{1,V_0\}^{\frac{1}{p}}\|u\|_{\mathcal {W}_\varepsilon}.
\end{equation}
Together with the continuous embedding $\mathcal {W}_{\varepsilon_n}\hookrightarrow W^{s,p}(\mathbb{R}^N),$
we get that $\{\mathfrak u_n\}_n$ is bounded in $\mathcal {W}.$

Next, we claim that there exist a sequence $\{{\tilde y}_n\}_n\subset{\mathbb{R}^N}$ and constants $R, \delta>0$ such that
\begin{align}\label{ctnewc4}
\liminf_{n \to  + \infty}{\int_{{B_R}\left( {{{\tilde y}_n}} \right)} {\left| {{u_n}} \right|} ^{q}}dx \geq \delta  > 0.
\end{align}
Suppose to the contrary, that for every $R>0,$ we deduce that
\[\mathop {\lim }\limits_{n \to  + \infty } \mathop {\sup }\limits_{y \in {\mathbb{R}^N}} {\int_{{B_R}\left( y \right)} {\left| {{\mathfrak u_n}} \right|} ^{q}}dx = 0.\]
Together with  Lemma \ref{lm3b}, we obtain that ${{\mathfrak u_n} \to 0}$ in $L^{\nu}(\mathbb{R}^{N}),$ for every $\nu \in \left( {p, + \infty } \right).$ Using Lemma \ref{lm1} and (\ref{ctnew3a}), we get
$$\lim_{n\to\infty}\int\limits_{\mathbb R^N}\Big[\frac{1}{|x|^{\mu}}*G(\varepsilon y, \mathfrak u_n(y))\Big] G(\varepsilon x,\mathfrak u_n(x))dx=0.$$
Together with $\mathfrak u_n\in\mathscr{N}_{\varepsilon_n},$ we have $\mathfrak u_n
\to 0$ in $\mathcal {W}_{\varepsilon}.$ Hence $\mathcal {I}_{\varepsilon}(u_n)
\to 0$, which is impossible due to $c_{Z_0}>0.$ We now suppose ${v_n}\left( x \right) = {\mathfrak u_n}\left( {x + {{\tilde y}_n}} \right),$ so ${\left\{ {{v_n}} \right\}_{n}}$  is bounded in $\mathcal {W}_{\varepsilon}.$  Therefore, we may assume that $v_{n}\rightharpoonup v$ in $\mathcal {W},$ as $n\rightarrow\infty.$ It follows from \eqref{ctnewc4} that $v\ne 0.$

Let $t_n>0$ be such that $\tilde v_n:=t_n v_n\in \mathscr{N}_{Z_0}$ and
$y_n:=\varepsilon_n \tilde y_n$. Applying Lemma \ref{lm8.1}, for every $n$ there exists a unique $t_{\mathfrak u_n}> 0$ such that $\mathcal {I}_{\varepsilon_n}(t_{\mathfrak u_n}\mathfrak u_n)=\text{sup}_{t\ge 0}\mathcal {I}_{\varepsilon_n}(t\mathfrak u_n).$ Then $t_{\mathfrak u_n}\mathfrak u_n\in \mathscr{N}_{\varepsilon_n},$ which yields
$t_{\mathfrak u_n}=1,$ due to $\mathfrak u_n\in\mathscr{N}_{\varepsilon_n}$. Therefore, $\sup_{t\ge 0}\mathcal {I}_{\varepsilon_n}(t\mathfrak u_n)=\mathcal {I}_{\varepsilon_n}(\mathfrak u_n)$.
By the change of variable $z=x+\tilde y_n$, 
we deduce that
\begin{align*}
    c_{Z_0}\le \mathcal {E}_{Z_0}(\tilde v_n) &= \frac{1}{p}\|\tilde v_n\|_{Z_0,W^{s,p}(\mathbb{R}^N)}^p
    +\frac{1}{q}\|\tilde v_n\|_{Z_0,W^{s,q}(\mathbb{R}^N)}^q
   \quad -\int_{\R^N}\Big[\dfrac{1}{|x|^{\mu}}*F(\tilde v_n(y))\Big]F(\tilde v_n(x))dx\\
    &\le \frac{1}{p}\|\tilde v_n\|_{Z_0,W^{s,p}(\mathbb{R}^N)}^p
    +\frac{1}{q}\|\tilde v_n\|_{Z_0,W^{s,q}(\mathbb{R}^N)}^q \\
    & \quad -\int_{\R^N}\Big[\frac{1}{|x|^{\mu}}\ast
    G(\varepsilon_n y+y_n,\tilde v_n(y))\Big] G(\varepsilon_n x+y_n,v_n(x)) dx\\
    &=\mathcal {I}_{\varepsilon_n}(t_n u_n)\le 
    \mathcal {I}_{\varepsilon_n}(\mathfrak u_n)\le c_{Z_0}+o_n(1)
\end{align*}
which yields $\mathcal {E}_{Z_0}(\tilde v_n)\to c_{Z_0},$ as $n\to +\infty$. By  the fact that $\{\tilde v_n\}_n\subset \mathscr{N}_{Z_0}$ and $(f_3)$, we can pick $C_1>0$ such that $\left\|\tilde v_n\right\|_{Z_0}\le C_1,$ for every $n\in\N$. In addition, since $v_n\not\to 0$ in $\mathcal {W}$, there exists $\tilde{C}_1>0$ such that $\left\|v_n\right\|_{Z_0}\ge \tilde{C}_1>0,$ for every $n\in\N$. Therefore,
$$
\tilde{C}_1t_n\le \left\|t_nv_n\right\|_{Z_0,\mathcal {W}}=\left\|\tilde
v_n\right\|_{Z_0,\mathcal {W}}\le C_1
$$
which yields $t_n\le \frac{C_1}{\tilde{C}_1},$ for every $n\in\N$. Consequently, going to a subsequence if necessary, we suppose that $t_n\to t_0\ge 0$ and $\tilde v_n\rightharpoonup \tilde v:=t_0v\not\equiv 0$ in $\mathcal {W}$ and $\tilde v_n\to \tilde v$ a.e. in $\R^N$. If $t_0=0$, then $\tilde v_n\to 0$ in $\mathcal {W}.$ Therefore  $\mathcal {E}_{Z_0}(\tilde v)\to 0$, which is impossible since  $c_{Z_0}>0,$ so we get that  $t_0>0$. Arguing  as Proposition \ref{pro1}, we obtain $\mathcal {E}'_{Z_0}(\tilde v)=0$.

In the sequel, we shall prove that
\begin{equation}\label{nm40a}
    \lim_{n\to+\infty}\left\|\tilde v_n\right\|_{Z_0,\mathcal {W}}=\left\|\tilde v\right\|_{Z_0,\mathcal {W}}.
\end{equation}

Invoking the Fatou lemma, we can deduce
\begin{equation}\label{nm41a}
    \left\|\tilde v\right\|_{Z_0,\mathcal {W}}\le \liminf_{n\to\infty}\left\|\tilde v_n\right\|_{Z_0,\mathcal {W}}.
\end{equation}
Assume to the contrary, that
$$
\left\|\tilde v\right\|_{Z_0,\mathcal {W}}<\limsup_{n\to\infty}\left\|\tilde v_n\right\|_{Z_0,\mathcal {W}}.
$$
In such a case we would get
\begin{align*}
    c_{Z_0}+o_n(1)&=\mathcal {E}_{Z_0}(\tilde v_n)-\frac{1}{2q}\langle \mathcal {E}'_{Z_0}(\tilde v_n), \tilde v_n\rangle\\
    &=(\frac{1}{p}-\frac{1}{2q})\left\|\tilde v_n\right\|_{Z_0,W^{s,q}(\mathbb{R}^N)}^q
    +\frac{1}{2q}\left\|\tilde v_n\right\|_{Z_0,W^{s,q}(\mathbb{R}^N)}^q \\
    & \quad
    +\frac{1}{2}\int_{\R^N}\Big[\frac{1}{|x|^{\mu}}*F(\tilde v_n(y))\Big]\left(\frac{1}{q}f(\tilde v_n)\tilde v_n-F(\tilde v_n)\right)dx
\end{align*}
and, by $(f_3)$ and the Fatou lemma, we would have
\begin{align*}
    c_{Z_0}&\ge \limsup_{n\to+\infty}\Big[(\frac{1}{p}-\frac{1}{2q})\left\|\tilde v_n\right\|_{Z_0,W^{s,q}(\mathbb{R}^N)}^q
    +\frac{1}{2q}\left\|\tilde v_n\right\|_{Z_0,W^{s,q}(\mathbb{R}^N)}^q\Big] \\
    & \quad +\frac{1}{2}\liminf_{n\to +\infty}\int_{\R^N}\Big[\frac{1}{|x|^{\mu}}*F(\tilde v_n(y))\Big]\left(\frac{1}{q}f(\tilde v_n)\tilde v_n-F(\tilde v_n)\right)dx\\
    &>(\frac{1}{p}-\frac{1}{2q})\left\|\tilde v\right\|_{Z_0,W^{s,q}(\mathbb{R}^N)}^q
    +\frac{1}{2q}\left\|\tilde v\right\|_{Z_0,W^{s,q}(\mathbb{R}^N)}^q \\
    & \quad + \frac{1}{2}\int_{\R^N}\Big[\frac{1}{|x|^{\mu}}*F(\tilde v(y))\Big]\left(\frac{1}{q}f(\tilde v)\tilde v-F(\tilde v)\right)dx \\
    &=\mathcal {E}_{Z_0}(\tilde v)-\frac{1}{2q}\mathcal {E}'_{Z_0}(\tilde v)(\tilde v)
    =\mathcal {E}_{Z_0}(\tilde v)\ge c_{Z_0}
\end{align*}
which is a contradiction. Therefore, $w_n\rightharpoonup w$ in $\mathcal {W}$ and \eqref{nm40a} implies $\tilde v_n\to \tilde v$ in $\mathcal {W}.$  Moreover, $v_n\to v$ in $\mathcal {W},$ as $n\to +\infty$.

In order to complete  the proof of this lemma, we explore ${y_n} = {\varepsilon _n}{y_n}.$ We claim that $\{y_n\}_n$ allows a subsequence, still denoted the same, such that ${{y_n} \to {y_0}},$ for some ${{y_0} \in \mathscr{M}.}$  In the sequel, we have to verify that the following two claims hold.

 \smallskip

\textbf{Claim 1.} $\{y_n\}_n$ is bounded.

We shall argue by  contradiction. Assume that, up to a subsequence, $\left| {{y_n}} \right| \to \infty, $ as $n\rightarrow\infty.$ Since $\left\langle {{\mathcal {I}_{{\varepsilon _n}}'}(\mathfrak u_n),{\mathfrak u_n}} \right\rangle  = 0$ and ${\mathcal {I}_{{\varepsilon _n}}}(\mathfrak u_n)\rightarrow{c_{{V_0}}},$ by Lemma \ref{lma}, we can infer that there exists $\hat{C_{0}}\in(0,\frac{\hbar_0}{2})$ such that
\[ {\left| {\frac{1}{{{{\left| x \right|}^\mu }}}*G\left( {\varepsilon y,\mathfrak u_n} \right)} \right| } < \hat{{C_0}}.\]
Fixed $R>0$ such that $\Lambda  \subset {B_R}\left( 0 \right),$ and assume that $\left| {{y_n}} \right| > 2R.$ Therefore,
\begin{align}\label{mo1}
\left| {{\varepsilon _n}x + {y_n}} \right| \geq \left| {{y_n}} \right| - \left| {{\varepsilon _n}x} \right| > R\
\hbox{ for every }
x \in {B_{\frac{R}{{{\varepsilon _n}}}}}\left( 0 \right).
\end{align}
Note that $\mathfrak u_n\in\mathscr{N}_{\varepsilon_n},$  so we have
\begin{align*}
\|\mathfrak u_n\|_{Z_0,W^{s,p}(\mathbb R^N)}^{p}
+\|\mathfrak u_n\|_{Z_0,W^{s,q}(\mathbb R^N)}^{q}
&  \le \|\mathfrak u_n\|_{W^{s,p}_{Z,\varepsilon_n}}^{p}+ \|\mathfrak u_n\|_{W^{s,q}_{Z,\varepsilon_n}}^{q}
 =\int\limits_{\mathbb R^N}\Big[\frac{1}{|x|^{\mu}}*G(\varepsilon y,\mathfrak u_n)\Big]g(\varepsilon_nx,\mathfrak u_n)\mathfrak u_ndx.
\end{align*}
Using the change of variables $x\mapsto x+\tilde y_n$ and $y\mapsto y+\tilde y_n,$ we get
\begin{align*}
\|v_n\|_{Z_0,W^{s,p}(\mathbb R^N)}^{p}
+\|v_n\|_{Z_0,W^{s,q}(\mathbb R^N)}^{q}
& =\|\mathfrak u_n\|_{W^{s,p}_{Z,\varepsilon_n}(\mathbb{R}^N)}^{p}+ \|\mathfrak u_n\|_{W^{s,q}_{Z,\varepsilon_n}(\mathbb{R}^N)}^{q} \\
&\le \int\limits_{\mathbb R^N}\Big[\frac{1}{|x|^{\mu}}*G(\varepsilon y,\mathfrak u_n)\Big]g(\varepsilon_nx,\mathfrak u_n)\mathfrak u_ndx\\
&=\int\limits_{\mathbb R^N}\Big[\frac{1}{|x|^{\mu}}*G(\varepsilon y+y_n,v_n)\Big]g(\varepsilon_nx+y_n,v_n)v_ndx.
\end{align*}
By (\ref{mo1}), the definition of $g,$ $f\left( t \right) \leq \frac{{{Z_0}}}{\hbar_0}{t^{p - 1}},$ $\hat{{C_0}} \in \left( {0,\frac{\hbar_0}{2}} \right),{v_n} \to v$ in $\mathcal {W}_{\varepsilon},$ the Dominated Convergence Theorem, we have that
\begin{align*}
\|v_n\|_{Z_0,W^{s,p}(\mathbb R^N)}^{p}
+\|v_n\|_{Z_0,W^{s,q}(\mathbb R^N)}^{q}
&\nonumber\leq \hat{{C_0}}\int_{{\mathbb{R}^N}} {g\left( {{\varepsilon _n}x + {y_n},{v_n}} \right)} {v_n}dx\\
&\nonumber \leq \hat{{C_0}}\int_{B_{\frac{R}{\epsilon_{n}}}(0)} {\hat{f}\left( {{v_n}} \right)} {v_n}dx + \hat{{C_0}}\int_{B_{\frac{R}{\epsilon_{n}}}^{c}(0)}f\left({v_n} \right) {v_n}dx\\
&\nonumber\leq\frac{1}{2}\int_{B_{\frac{R}{\epsilon_{n}}}\left( 0
\right)} {{Z_0}} ({\left| {{v_n}} \right|^p}+|v_n|^q)dx +
{o_n}\left( 1 \right)
\end{align*}
which gives $$\|v_n\|_{Z_0,W^{s,p}(\mathbb R^N)}^{p}
+\|v_n\|_{Z_0,W^{s,q}(\mathbb R^N)}^{q} = {o_n}\left( 1 \right).$$  Therefore, we have that $\{y_n\}_n$ is bounded in $\mathbb{R}^N$.

 \smallskip

\textbf{Claim 2.}  $y_0\in \mathscr{M}.$

By Claim 1, up to a
subsequence, we can suppose that ${{y_n} \to {y_0}.}$ Once  ${{y_0} \notin\bar{\Omega}}$, which implies
 the closure of $\Omega$, we can argue as above to get $v_{n}\rightarrow 0$ in $\mathcal {W}_{\varepsilon},$ which is impossible. Therefore, we obtain ${{y_0} \in \bar{\Omega} }.$ Now, suppose by contradiction that $Z\left( {{y_0}} \right) > {Z_0},$ then by using ${\tilde{v}_n} \to v$ in $\mathcal {W}$ and the Fatou lemma, we can deduce that
\begin{align}\label{o7}
{c_{{Z_0}}}\nonumber = {\mathcal {E}_{{Z_0}}}\left( \tilde{v }\right)
&< \liminf\Big[
\frac{1}{p}\|\tilde{v}_n\|_{Z_0,W^{s,p}(\mathbb R^N)}^{p}
+\frac{1}{q}\|\tilde{v}_n\|_{Z_0,W^{s,q}(\mathbb R^N)}^{q} \\
&\nonumber\quad - \frac{1}{2}\int_{{\mathbb{R}^N}} {\left( {\frac{1}{{{{\left| x \right|}^\mu }}}\ast F\left( {{\tilde{v}_n}} \right)} \right)} F\left( {{\tilde{v}_n}} \right)dx\Big]\\
&\nonumber\leq\liminf_{n \to \infty }{\mathcal {I}_{{\varepsilon _n}}}\left( {t_n{\mathfrak u_n}} \right)\leq\nonumber\liminf_{n \to \infty }{\mathcal {I}_{{\varepsilon _n}}}\left(\mathfrak u_n \right)=c_{Z_0}
\end{align}
which is impossible. Therefore, $Z\left( {{y_0}} \right) = {Z_0}$ and ${{y_0} \in \bar{\Omega} .}$
 Thanks to $(\mathcal {Z}_2)$, $y_0\notin \partial\Omega$, and thus
$y_0\in \mathscr{M}.$ This completes the proof of Lemma \ref{lm10}.
\end{proof}

Let $h(\varepsilon)$ be any positive function satisfying $h(\varepsilon)\rightarrow 0,$ as $\varepsilon\rightarrow 0.$ Define
\[\tilde{{\mathscr{N}}}_\varepsilon  = \left\{ {u \in {\mathscr{N}_\varepsilon }:{\mathcal {I}_\varepsilon }(u) \leq {c_{{Z_0}}} + h(\varepsilon)} \right\}.\]
For any $y \in \mathscr{M},$ we deduce from Lemma \ref{lm8} that $h\left( \varepsilon  \right) = \left| {{\mathcal {I}_\varepsilon }\left( {{\Phi _\varepsilon }\left( y \right)} \right) - {c_{{Z_0}}}} \right| \to 0,$ as $\varepsilon\rightarrow 0.$ Thus, ${{\Phi _\varepsilon }\left( y \right) \in {\tilde{\mathscr{N}_\varepsilon }}}$ and ${\tilde{{\mathscr{N}_\varepsilon }} \ne \phi}$ for every $\varepsilon>0.$

 \begin{lemma}\label{lm11}(see Thin \cite[Lemma 16]{T22})
 For every $\delta>0,$
 $$ \lim\limits_{\varepsilon\to 0^{+}}\sup\limits_{u\in
 \tilde {\mathscr{N}}_{\varepsilon}}
 \mbox{\rm dist}(\beta_{\varepsilon}(u),\mathscr{M}_{\delta})=0.$$
 \end{lemma}

\begin{lemma}\label{lm12}
Suppose that conditions $(\mathcal {Z}_1)-(\mathcal {Z}_2)$ and $(f_1)-(f_5)$ hold and
denote by $\mathfrak v_n$  a
 nontrivial nonnegative solution in~$\mathbb R^N$ of
 \begin{align}\label{nm45a}
 (-\Delta)_{N/s}^{s}\mathfrak v_n+ (-\Delta)_{q}^{s}\mathfrak v_n
 +Z_n(x)(|\mathfrak v_n|^{\frac{N}{s}-2}\mathfrak v_n+|\mathfrak v_n|^{q-2}\mathfrak v_n)  =[|x|^{-\mu}
 \ast  F(\mathfrak v_n)]g(\varepsilon_nx+\varepsilon_n\tilde{y}_n,\mathfrak v_n),
 \end{align}
 where $Z_n(x)=Z(\varepsilon_n x+\varepsilon_n \tilde y_n)$ and $\varepsilon_n\tilde y_n\to y\in\mathscr{M}.$ Then, if $(\mathfrak v_n)_n$ is a bounded sequence in $\mathcal {W}$ satisfying
$$\limsup_{n\rightarrow\infty}\|u_n\|_{\vartheta, W^{s,p}(\mathbb{R}^N)}^{N/(N-s)}<
\frac{\beta_\ast\mathfrak d_\ast^{s/(N-s)}}{\mathfrak c\alpha_0}, \ \ \mbox{with}
 \ \ 0<\beta_\ast<\alpha_\ast$$
for a suitable constant   $\mathfrak c>1$
 and if $\mathfrak v_n\to \mathfrak v$ in $\mathcal {W}$, then each $\mathfrak v_n\in L^{\infty}(\mathbb R^N)$ and there
 exists $C>0$ such that $\|\mathfrak v_n\|_{L^{\infty}(\mathbb R^N)}\le C$
 for every $n.$ Moreover,
 $$ \lim\limits_{|x|\to\infty}\mathfrak v_n(x)=0,\quad \text{uniformly in}\; n.$$
 \end{lemma}

 \begin{proof}
In view of $\mathcal {I}_{\varepsilon_n}(u_n)\leq c_{Z_0}+ h(\varepsilon)$, with $ h(\varepsilon)\rightarrow0,$
as $n\rightarrow\infty.$ We argue as in the proof
of Lemma \ref{lm11} to show that $\mathcal {I}(\varepsilon_n)(u_n)\rightarrow c_{Z_0}.$
Then
by  Lemma \ref{lm10}, there exists $\{\tilde{y}_n\}\subset \mathbb{R}^N$
such that $v_n=u_n(\cdot+\tilde{y}_n)$ strongly converges in $\mathcal {W}$
and $\varepsilon_n\tilde{y}_n\rightarrow y_0\in \mathscr{M}.$
By the boundedness of $\{v_n\}_n$
in $\mathcal {W},$ we can proceed as in the proof of Lemma \ref{lma} to obtain that there exists
$\hat{C}_0>0$ such that
$$\frac{1}{|x|^{\mu}}\ast G(\varepsilon_nx+\varepsilon_n\tilde{y},v_n)\leq \hat{C}_0.$$
Repeating the same Moser iteration argument developed in the proof
of   Liang et al.~\cite[Lemma 5.5]{Liang1}, we have that
$\|\mathfrak v_n\|_{L^{\infty}(\mathbb R^N)}\le C$ for every $n\in
\mathbb{N}$. Now, we note that $v_n$ satisfies problem
\eqref{nm45a}.

Using  Ambrosio and R$\breve{\mbox{a}}$dulescu
\cite[Corollary 2.1 ]{Amb3}
 and the fact that $\mathfrak v_n$ is uniformly bounded in $L^{\infty}(\mathbb{R}^N)\cap \mathcal {W}$, we
can conclude that $\mathfrak v_n(x)\rightarrow 0$ as $|x|\rightarrow
\infty$ uniformly in $n\in\mathbb{N}.$
This completes the proof of Lemma \ref{lm12}.
 \end{proof}

\begin{proof}[Proof of Theorem \ref{th2}]
Using the similar arguments to the proof of
Ambrosio \cite[Theorem 5.2]{DEF}
and~\cite[Theorem 1.1]{PPP}. We define ${\alpha _\varepsilon }:\mathscr{M} \to {\mathbb{S}_\varepsilon }$ by setting ${\alpha _\varepsilon }\left( y \right) = m_\varepsilon ^{ - 1}\left( {{\Phi _\varepsilon }\left( y \right)} \right)$ for every $\varepsilon  > 0.$ Applying Lemma \ref{lm8} and the definition of ${\Phi _\varepsilon },$  we obtain that
\begin{center}
$\mathop {\lim }\limits_{\varepsilon  \to 0} {\psi_\varepsilon }\left( {{\alpha _\varepsilon }\left( y \right)} \right) = \mathop {\lim }\limits_{\varepsilon  \to 0} {\mathcal {I}_\varepsilon }\left( {{\Phi _\varepsilon }\left( y \right)} \right) = c_{Z_0},$ uniformly in $y \in \mathscr{M}.$
\end{center}
Therefore, there exists $\tilde{\varepsilon}>0$ such that
${\bar{\mathbb{S}}_\varepsilon }: = \left\{ {w \in {\mathbb{S}_\varepsilon }:{\psi_\varepsilon }\left( w \right) \leq c_{Z_0} + h\left( \varepsilon  \right)} \right\} \ne \phi, $ for every $\varepsilon  \in \left( {0,\tilde{\varepsilon}} \right).$ With the aid of Lemma \ref{lm8.1}-$(iii)$, Lemma \ref{lm8}
and Lemma \ref{lm11}, for every $\delta>0,$ there exists ${\tilde{\varepsilon } = {{\tilde{\varepsilon} }_\delta } > 0}$ such that the diagram of continuous mappings
$$\mathscr{M}\mathop  \to \limits^{{\Phi _\varepsilon }} {{\tilde {\mathscr{N}}}_\varepsilon }\mathop  \to \limits^{m_\varepsilon ^{ - 1}} {\mathbb{S}_\varepsilon }\mathop  \to \limits^{{m_\varepsilon }} {{\tilde {\mathscr{N}}}_\varepsilon }\mathop  \to \limits^{{\beta _\varepsilon }} {\mathscr{M}_\delta } = {{ \tilde{\varepsilon} }_\delta } > 0$$
is well-defined, for every $\varepsilon  \in \left( {0,\tilde{\varepsilon} } \right).$
Invoke Lemma \ref{lm9} and take a function $\varpi\left( {\varepsilon ,y} \right)$ with $\left| {\varpi\left( {\varepsilon ,y} \right)} \right| < \frac{\delta }{2}$ uniformly in $y \in \mathscr{M},$ for every $\varepsilon  \in \left( {0,\tilde{\varepsilon} } \right)$ such that ${\beta _\varepsilon }\left( {{\Phi _\varepsilon }\left( y \right)} \right) = y +\varpi\left( {\varepsilon ,y} \right),$ for every $y \in \mathscr{M}.$ Therefore, we obtain that $\mathcal {F}\left( {t,y} \right) = y + \left( {1 - t} \right)\varpi\left( {\varepsilon ,y} \right)$ with $\left( {t,y} \right) \in \left[ {0,1} \right] \times\mathscr{M}$ is a homotopy between ${\beta _\varepsilon } \circ {\Phi _\varepsilon } = \left( {{\beta _\varepsilon } \circ {m_\varepsilon }} \right) \circ {\alpha _\varepsilon }$ and the inclusion map id:$\mathscr{M}\to {\mathscr{M}_\delta }.$ Together with~\cite[Lemma 6.3.21]{Amb1}, we obtain that
\begin{align}\label{3.12}
ca{t_{{\tilde{\mathbb{S}}_\varepsilon }}}\left( {\tilde{{\mathbb{S}}_\varepsilon }} \right) \geq ca{t_{{\mathscr{M}_\delta }}}\left(\mathscr{M}\right).
\end{align}

In what follows, we  choose a function $h\left( \varepsilon  \right) > 0$ such that $h\left( \varepsilon  \right) \to 0,$ as $\varepsilon  \to 0$ and such that ${c_{Z_0} + h\left( \varepsilon  \right)}$ is not a critical level for ${\mathcal {I}_\varepsilon }.$ From Lemma \ref{lm3.7}, we see that
${{\mathcal {I}_\varepsilon }}$ satisfies the Palais-Smale condition in ${{{\tilde{\mathbb{S}}}_\varepsilon }}$ as $\varepsilon >0.$ Invoking
Ambrosio \cite[Theorem 6.3.20]{Amb1},
we get that ${{\psi_\varepsilon }}$ has at least
$ca{t_{{{\tilde{\mathbb{S}}}_\varepsilon }}}(\tilde{\mathbb{S}}_\varepsilon)$ critical points on $\tilde{\mathbb{S}}_\varepsilon.$ Consequently, by Lemma \ref{lm8.2} and (\ref{3.12}), we deduce
 that ${\mathcal {I}_\varepsilon }$ has at least $ca{t_{{\mathscr{M}_\delta }}}\left(\mathscr{M}\right)$ critical points.

Let $u_{\varepsilon_n}$ be a solution of problem
$(\mathcal {Q}_{\varepsilon_n}),$ then $\mathfrak v_n(x)=u_{\varepsilon_n}(x+\tilde y_n)$ is also a solution of problem \eqref{nm45a}.  Moreover, there
 exists $\mathfrak v\in \mathcal {W}$,
  such that, up to a subsequence,
 $\mathfrak v_n\to \mathfrak v$  in $\mathcal {W}$
 and $y_n=\varepsilon_n\tilde y_n\to y\in \mathscr{M}$ by Lemma~\ref{lm10}.

 We claim that there exists $\bar{\delta}>0$
 such that $\|v_n\|_{L^{\infty}(\mathbb R^N)}\ge\bar{\delta},$
 for every $n$ large enough. In fact, \eqref{ctnewc4} in the proof
  of Lemma~\ref{lm10} implies
  \begin{align}\label{54a}
  0<\frac{\beta}{2}\le \int\limits_{B_r(0)}|\mathfrak v_n|^{N/s}dx\le |B_r(0)|\|\mathfrak v_n\|_{L^{\infty}(\mathbb R^N)}^{N/s}
  \end{align}
for every $n$ large enough. Therefore,  we choose $\bar{\delta}=\Big(\frac{\beta}{2|B_r(0)|}\Big)^{s/N}.$
Applying the fact that $\mathfrak v_n\to \mathfrak v$ in $\mathcal {W},$ we have
 $\lim\limits_{|x|\to\infty}
 \mathfrak v_n(x)=0,$ uniformly in $n$ by Lemma~\ref{lm12}.

 Let $\mathfrak q_n$ be a global maximum point of $\mathfrak v_n.$ Invoking Lemma~\ref{lm12}, we see that there exists $R>0$ such that $|\mathfrak q_n|\le R,$ for every $n$.
 Consequently, the maximum point of $u_{\varepsilon_n}$ is denoted  by $\mathfrak z_{\varepsilon_n}=\mathfrak q_n+\tilde y_n.$ Furthermore, problem \eqref{pr2} possesses  a nontrivial nonnegative solution $w_{\varepsilon}(x)=u_{\varepsilon}(x/{\varepsilon}).$ Therefore, the maximum points $\zeta_{\varepsilon}$ of $w_{\varepsilon}$ and $\mathfrak z_{\varepsilon}$
  of $u_{\varepsilon}$ satisfy
 $\zeta_{\varepsilon}=\varepsilon \mathfrak z_{\varepsilon}.$ We see that
 $$ \lim_{\varepsilon\to0^{+}}Z(\zeta_{\varepsilon})
 =\lim_{n\to\infty}Z(\varepsilon_n\mathfrak z_{\varepsilon_n})=Z_0.$$%
 This completes the proof of  Theorem \ref{th2}.
\end{proof}

\begin{proof}[Proof of Theorem $\ref{th2a}$]
We know that $w_{\varepsilon}(x)=u_{\varepsilon}(x/{\varepsilon})$ is a nontrivial nonnegative solution of problem \eqref{pr2}.
Set
$$\mathfrak v_{\varepsilon_n}=w_{\varepsilon_n}(\varepsilon_n
\cdot+\eta_{\varepsilon_n})=u_{\varepsilon_n}(\cdot+\mathfrak z_{\varepsilon_n}).$$
Therefore, Lemma \ref{lm10} yields that $(\mathfrak v_{\varepsilon_n})_n\rightarrow \mathfrak v$
in $\mathcal {W}$ and $\mathfrak v$ is a ground state solution of
the following equation
$$(-\Delta)_p^su+(-\Delta)_q^s+Z_{0}(|u|^{p-2}v+|u|^{q-2}u)
=\Big[{|x|^{-\mu}}*F(u(y))\Big]f(u)\;\text{in}\;\mathbb R^N.$$ This
completes the proof of Theorem \ref{th2a}.
\end{proof}

\begin{proof}[Proofs of Corollaries \ref{cornew10} and \ref{cornew1}]
Apply a similar discussion as in Liang et al.~\cite{Liang1}.
\end{proof} 
\hfill 

\noindent{\bf Acknowledgements.}
The first two authors were supported  by the Science and Technology
Development Plan Project of Jilin Province, China (No.
20230101287JC), the Research Foundation of Department
of Education of Jilin Province ( No. JJKH20251034KJ), the Young  outstanding talents project of Scientific Innovation and
entrepreneurship in Jilin (No. 20240601048RC). The third author was supported by the Slovenian Research and Innovation Agency program P1-0292 and grants J1-4031, J1-4001, N1-0278, N1-0114, and N1-0083.
We thank the referee for comments and suggestions.
\\


\begin{thebibliography}{x}

\bibitem{ACTY} \textsc{C.O. Alves, D. Cassani, C. Tarsi, M. Yang},
Existence and concentration of ground state solutions for a critical nonlocal Schr\"odinger equation in $\mathbb R^2,$
J. Differential Equations 261 (2016) 1933--1972.

\bibitem{PPP} \textsc{V. Ambrosio,}
On the multiplicity and concentration of positive solutions for a $p$-fractional
Choquard equation in $\mathbb{R}^{N}$,
Comput. Math. Appl. 78 (2019) 2593--2617.

\bibitem{DEF} \textsc{V. Ambrosio,}
Multiplicity and concentration results for a fractional Choquard equation via penalization method,
Potential Anal. 50 (2019) 55--82.

\bibitem{Amb1}  \textsc{V. Ambrosio},
Nonlinear fractional Schr\"{o}dinger equations in $\mathbb{R}^N$,
Frontiers in Elliptic and Parabolic Problems. Birkh\"{a}user/Springer, Cham, 2021.

\bibitem{A0} \textsc{V. Ambrosio},
Multiple concentrating solutions for a fractional $(p, q)$-Choquard equation,
Adv. Nonlinear Stud. 24 (2024) 510--541.

\bibitem{Amb}  \textsc{V. Ambrosio, T. Isernia,}
Multiplicity and concentration results for some nonlinear Schr\"{o}dinger equations with the fractional $p$-Laplacian,
Discrete Contin. Dyn. Syst. 38 (2018) 5835--5881.

\bibitem{Amb3}   \textsc{V. Ambrosio, V.D. R\u{a}dulescu},
Fractional double-phase patterns: concentration and multiplicity of solutions,
J. Math. Pures Appl. 142 (2020) 101--145.

 \bibitem{AS} \textsc{ S.N. Antontsev, S.I. Shmarev},
 Elliptic equations and systems with nonstandard growth conditions: Existence, uniqueness and localization properties of solutions,
Nonlinear Anal. 65 (2006) 722--755.

\bibitem{BDaFP} \textsc{V. Benci, P. D'Avenia, D. Fortunato, L. Pisani},
Solitons in several space dimensions: Derrick's problem and infinitely many solutions,
Arch. Ration. Mech. Anal.  154 (2000)  297--324.

\bibitem{BAP} \textsc{D. Bonheure, P. d'Avenia,  A. Pomponio},
On the electrostatic Born-Infeld equation with extended charges,
Commun. Math. Phys. 346 (2016)  877--906.

\bibitem{BL} \textsc{M. Born,  L. Infeld},
Foundations of the new field theory,
Nature 132 (1933) 1004.

\bibitem{BL1} \textsc{H. Br\'{e}zis, E. Lieb},
A relation between pointwise convergence of functions and convergence of functionals,
Proc. Amer. Math. Soc. 88 (1983) 486--490.

\bibitem{CLY} \textsc{S. Chen, L. Li, Z. Yang,}
Multiplicity and concentration of nontrivial nonnegative solutions for a fractional Choquard equation with critical exponent,
Rev. R. Acad. Cienc. Exactas Fs. Nat. Ser. A Mat. RACSAM 114(2020) Paper No. 33, 35 pp.

\bibitem{CTS22} \textsc{S. Chen, M. Shu, X. Tang, L. Wen,}
Planar Schr\"odinger-Poisson system with critical exponential growth in the zero mass case,
J. Differential Equations 327 (2022) 448--480.

\bibitem{chen1} \textsc{Y. Chen,  Z. Yang},
The existence of multiple solutions for a class of upper critical Choquard equation in a bounded domain,
Demonstr. Math. 57 (2024) 20230152.

\bibitem{CV} \textsc{L. Cherfil,  V. Il'yasov},
On the stationary solutions of generalized reaction diffusion equations with $p\&q$-Laplacian,
Commun. Pure Appl. Anal. 1 (2004) 1--14.

\bibitem{Cin} \textsc{S. Cingolani,  K. Tanaka},
Semi-classical states for the nonlinear Choquard equations: existence, multiplicity and concentration at a potential well,
Rev. Mat. Iberoam. 35 (2019) 1885--1924.

\bibitem{CAB} \textsc{R. Clemente, J.C.D. Albuquerque, E. Barboza},
Existence of solutions for a fractional Choquard-type equation in $\mathbb R^N$ with critical exponential growth,
Z. Angew. Math. Phys. 72 (2021) 16 pp.

\bibitem{EM21} \textsc{E. de  S. B\"oer, O.H. Miyagaki,}
Existence and multiplicity of solutions for the fractional $p$-Laplacian Choquard logarithmic equation involving a nonlinearity with exponential critical and subcritical growth,
J. Math. Phys. 62 (2021) 051507.

\bibitem{Pi1} \textsc{M. del Pino, P. Felmer},
Local mountain passes for semilinear elliptic problems in unbounded domains,
Calc. Var. Partial Differential Equations  4 (1996) 121--137.

 \bibitem{Ne1} \textsc{E. Di Nezza, G.Palatucci, E.Valdinoci},
 Hitchhiker's guide to the fractional Sobolev spaces,
 Bull. Sci. Math. 136 (2012) 521--573.

\bibitem{FW86} \textsc{A. Floer, A. Weinstein},
Non spreading wave packets for
the cubic Schr\"odinger equation with a bounded potential,
J. Funct. Anal. 69 (1986) 397--408.

\bibitem{LY} \textsc{Q.  Li, Z. Yang,}
Multiple solutions for a class of fractional quasi-linear equations with critical exponential growth in $\mathbb R^N,$
Complex Var. Elliptic Equ. 61 (2016) 969--983.

\bibitem{Liang1} \textsc{S. Liang, P. Pucci, T. Van Nguyen},
Multiplicity and concentration results for some fractional double phase Choquard equation with exponential growth, 
Asymptot. Anal. (2025), publ. online. DOI: 10.1177/09217134251319160

\bibitem{SLN23} \textsc{S. Liang, S. Shi,  T.V. Nguyen},
Multiplicity and concentration properties for fractional Choquard equations with exponential growth,
J. Geom. Anal.  34 (2024) 367.

\bibitem{L76} \textsc{E.H. Lieb,}
Existence and uniqueness of the minimizing solution of Choquard's nonlinear equation,
Stud. Appl. Math. 57  (1976/77)  93--105.

\bibitem{Lie} \textsc{E. Lieb, M. Loss,}
 {\em Analysis},  Grad. Stud. Math. 14, American Mathematical Society, Providence, RI, 2001.

\bibitem{PL80} \textsc{P.L. Lions,} The Choquard equation and related questions,
Nonlinear Anal.: Theory, Methods Appl. 4 (1980) 1063--1072.

\bibitem{BTR23} \textsc{G. Molica  Bisci, N.V. Thin, V.D. R\u adulescu},
Concentration phenomena for  fractional double phase equations with Choquard reaction, preprint.

\bibitem{NVT1} \textsc{G. Molica  Bisci, N.V. Thin, L. Vilasi},
On a class of nonlocal Schr\"odinger equations with exponential growth,
Adv. Differential Equations 27 (2022) 571--610.

\bibitem{PRR}  \textsc{N.S. Papageorgiou, V.D. R\u{a}dulescu, D.D. Repov\v{s},}
Nonlinear Analysis - Theory and Applications,
Springer, Cham, 2019.

\bibitem{PR18} \textsc{E. Parini, B. Ruf,}
On the Moser-Trudinger inequality in fractional Sobolev-Slobodeckij spaces,
Atti Accad. Naz. Lincei Rend. Lincei Mat. Appl. 29 (2018) 315--319.

 \bibitem{PXZ} \textsc{P. Pucci, M. Xiang, B. Zhang,}
 Multiple solutions for nonhomogeneous Schr\"odinger-Kirchhoff type equations involving the fractional $p$-Laplacian in $\mathbb R^N$,
Calc. Var. Partial Differential Equations  54 (2015) 2785--2806.

\bibitem{LQ} \textsc{A. Szulkin, T. Weth},
The method of Nehari manifold, Handbook of Nonconvex Analysis and Applications,
edited by D. Y. Gao and D. Motreanu, International Press, Boston, 2010, pp. 597--632.

\bibitem{T22} \textsc{N.V. Thin},
Multiplicity and concentration of solutions to a fractional $p$-Laplace problem with exponential growth,
Ann. Fenn. Math. 47 (2022) 603--639.

\bibitem{TTL23} \textsc{N.V. Thin, P.T. Thuy. T.T.D. Linh},
Existence of solution for the $(p,q)$-fractional Laplacian equation with nonlocal Choquard reaction and exponential growth,
 Complex Var. Elliptic Equ. 69 (2024) 1949--1972.

\bibitem{NW96}\textsc{M. Willem},
Minimax Theorems,
Progress in Nonlinear Differential Equations and Their Applications, Birkh\"auser 1996.

\bibitem{Yang} \textsc{Z. Yang, F. Zhao},
Multiplicity and concentration behaviour of solutions for a fractional Choquard equation with critical growth,
Adv. Nonlinear Anal. 10 (2021) 732--774.

\bibitem{ZXN23}\textsc{B. Zhang, X. Han, N. Thin,}
Schr\"odinger-Kirchhoff-type problems involving the fractional $p$-Laplacian with exponential growth,
Appl. Anal. 102 (2023) 1942--1974.

\bibitem{CZ} \textsc{C. Zhang},
Trudinger-Moser inequalities in Fractional Sobolev-Slobodeckij spaces and multiplicity of weak solutions to the Fractional-Laplacian equation,
Adv. Nonlinear Stud. 19 (2019) 197--217.

\bibitem{zl} \textsc{L. Zhang, Y. Liu, J.J. Nieto, G. Wang},
Nonexistence of solutions to fractional parabolic problem with general nonlinearities,
Rend. Circ. Mat. Palermo, II. Ser 73 (2024) 551--562.

 \bibitem{ZZR23} \textsc{W. Zhang, J. Zhang, V.D. R\u adulescu},
 Concentrating solutions for singularly perturbed double phase problems with nonlocal reaction,
J. Differential Equations 347 (2023)  56--103.

 \bibitem{Zhang2} \textsc{X. Zhang, X. Sun, S. Liang,  V.T. Nguyen,}
 Existence and concentration of solutions to a Choquard equation involving fractional $p$-Laplace via penalization method,
 J. Geom. Anal. 34 (2024) 90.

 \bibitem{zuo} \textsc{J. Zuo,  C. Liu, C. Vetro}, Normalized solutions to the fractional Schr\"odinger equation with potential, Mediterr. J. Math. 20 (2023) 216.

\end{thebibliography}
\end{document}